\documentclass[11pt]{article}

\usepackage{jheppub}
\usepackage{amsmath,amssymb,amsfonts,graphicx,slashed,amsthm,mathtools,upgreek,enumerate,tensor,subfig,bm}
\usepackage{wasysym}
\usepackage[dvipsnames]{xcolor}
\usepackage{arydshln,soul}
\usepackage{listings}
\usepackage{comment,hhline}
\usepackage[utf8]{inputenc}
\usepackage[titletoc]{appendix}
\usepackage{hyperref}

\makeatletter
\def\@fpheader{\relax}
\makeatother

\renewcommand{\footnoterule}{\vfill\kern -3pt \hrule width 0.4\columnwidth \kern 2.6pt} 

\definecolor{codegreen}{rgb}{0,0.6,0}
\definecolor{codegray}{rgb}{0.5,0.5,0.5}
\definecolor{codepurple}{rgb}{0.6,0.5,0.8}
\definecolor{codeback}{rgb}{0.95,0.95,0.95}

\lstdefinestyle{base}{
  basicstyle=\ttfamily\color{black},
  moredelim=**[is][\color{codepurple}]{!}{!},
}

\lstdefinestyle{pystyle}{
    backgroundcolor=\color{codeback},
    commentstyle=\color{codegreen},
    keywordstyle=\color{magenta},
    numberstyle=\tiny\color{codegray},
    stringstyle=\color{codepurple},
    basicstyle=\ttfamily\footnotesize,
    breakatwhitespace=false,
    breaklines=true,
    captionpos=t,
    keepspaces=true,
    numbers=left,
    numbersep=5pt,
    showspaces=false,
    showstringspaces=false,
    showtabs=false,
    tabsize=2
}
\DeclareMathOperator{\slope}{slope}
\DeclareMathOperator{\Vol}{Vol}
\DeclareMathOperator{\inj}{inj}
\DeclareMathOperator{\Kh}{Kh}
\DeclareMathOperator{\KFH}{KFH}
\newtheorem{theorem}{Theorem}

\makeatletter 
\def\fullwidthdisplay{\displayindent\z@ \displaywidth\columnwidth}
\edef\@tempa{\noexpand\fullwidthdisplay\the\everydisplay}
\everydisplay\expandafter{\@tempa}
\makeatother

\newcommand\bigcircle{\ocircle}

\newcommand\be{\begin{equation}}
\newcommand\ee{\end{equation}}
\newcommand\bea{\begin{eqnarray}}
\newcommand\eea{\end{eqnarray}}

\newcommand\eref[1]{(\ref{#1})}
\newcommand\bc{\begin{center}}
\newcommand\ec{\end{center}}
\renewcommand\comment[1]{}

\numberwithin{equation}{section} 
\allowdisplaybreaks

\newcommand{\MH}[1]{{\color{orange}{[\textbf{MH:} #1}]}}

\newenvironment{nohyphens}{%
  \hyphenpenalty=10000
  \exhyphenpenalty=10000
  \sloppy %
}{\par}

\title{\textrm{Illuminating new and known relations between knot invariants}
}

\author[a]{\textrm{Jessica Craven}}
\author[b]{\!\!, \textrm{Mark Hughes}}
\author[c]{\!\!, \textrm{Vishnu Jejjala}}
\author[d]{\!\!, \textrm{Arjun Kar}}

\affiliation[\,a]{Division of Physics, Mathematics, and Astronomy (PMA), California Institute of Technology,\\
Pasadena, CA 91125, USA}
\affiliation[\,b]{Department of Mathematics, Brigham Young University,\\
275 TMCB, Provo, UT 84602, USA}
\affiliation[\,c]{Mandelstam Institute for Theoretical Physics, School of Physics, NITheCS, and CoE-MaSS,\\
University of the Witwatersrand, 1 Jan Smuts Avenue, Johannesburg, WITS 2050, South Africa}
\affiliation[\,d]{Department of Physics and Astronomy, University of British Columbia,\\
6224 Agricultural Road, Vancouver, BC V6T 1Z1, Canada}

\emailAdd{jcraven@caltech.edu}
\emailAdd{hughes@mathematics.byu.edu}
\emailAdd{v.jejjala@wits.ac.za}
\emailAdd{arjunkar@phas.ubc.ca}

\abstract{\begin{nohyphens}
We automate the process of machine learning correlations between knot invariants.
For nearly 200,000 distinct sets of input knot invariants together with an output invariant, we attempt to learn the output invariant by training a neural network on the input invariants.
Correlation between invariants is measured by the accuracy of the neural network prediction, and bipartite or tripartite correlations are sequentially filtered from the input invariant sets so that experiments with larger input sets are checking for true multipartite correlation.
We rediscover several known relationships between polynomial, homological, and hyperbolic knot invariants, while also finding novel correlations which are not explained by known results in knot theory.
These unexplained correlations strengthen previous observations concerning links between Khovanov and knot Floer homology.
Our results also point to a new connection between quantum algebraic and hyperbolic invariants, similar to the generalized volume conjecture.
\end{nohyphens}}

\begin{document}

\maketitle
\parskip=.4\baselineskip

\section{Introduction}\label{sec:intro}
A knot $K$ is the image of an embedding $S^1 \hookrightarrow S^3$.
The same knot can, however, be drawn in $\mathbb{R}^2$ with over/undercrossing information in numerous ways, and these different knot diagrams are related to each other by enacting sequences of Reidemeister moves.
In general, it can be difficult to determine whether two knot diagrams represent the same knot because the required sequence of local moves may be complicated.
The ambiguity in representation is partially resolved by appealing to topological invariants, which are certain numbers, polynomials, or other algebraic structures that may be computed using a knot diagram but that are independent of the particular diagram chosen for a specific knot.
If two diagrams have different topological invariants, they represent different knots. 
The converse is not necessarily true.

The most common invariants have various mathematical and physical origins.
For example, the Jones polynomial is an algebraic invariant associated to the Hecke algebra of the braid group~\cite{jones85}.
It is equivalently described by skein relations and the Kauffman bracket~\cite{HKAUFFMAN1987395} or, as Witten showed, in quantum field theory as the unknot normalized vacuum expectation value of the Wilson loop operator in $SU(2)$ Chern--Simons gauge theory~\cite{Witten:1988hf}.
Evaluating the trace of the Wilson loop in representations of $SU(2)$ other than the fundamental one, we obtain colored Jones polynomials.
Promoting $SU(2)$ to $SU(N)$ generalizes the Jones polynomial to the HOMFLY-PT polynomial~\cite{Freyd:1985dx, przytycki2016invariants, Witten:1988hf}.
There are other polynomial invariants as well~\cite{Alexander:1923,alexander1928topological}.

Polynomial invariants in knot theory are typically Laurent polynomials with integer coefficients.
This fact about the Jones polynomial is explained by associating the coefficients to the dimensions of certain bigraded homology groups, thus realizing the Jones polynomial as the graded Euler characteristic of this homology theory. 
The bigrading in this Khovanov homology yields a new polynomial knot invariant, the Khovanov polynomial, a two variable polynomial whose powers are the homological and quantum gradings~\cite{khovanov2000, Bar_Natan_2002}.
Just as Khovanov homology categorifies the Jones polynomial, knot Floer homology categorifies the Alexander polynomial~\cite{ozsvath2004holomorphic}.
These homology theories lead to a new class of integer-valued topological invariants such as the Rasmussen $s$-invariant~\cite{rasmussen2010khovanov}, the Ozsv\'ath--Szab\'o $\tau$-invariant~\cite{ozsvath2004holomorphic}, and Hom's $\varepsilon$-invariant \cite{hom2014bordered}.
\footnote{
The HOMFLY-PT polynomial is categorified by a triply graded homology theory~\cite{khovanov2008matrix1, khovanov2008matrix2}.
New numerical invariants are likewise associated to this Khovanov--Rozansky homology~\cite{Gukov:2004hz,Dunfield:2005si,Gukov:2016gkn,Gukov:2017kmk,Gukov:2019mnk,Chandler:2022str}. }

Other invariants 
are associated only to hyperbolic knots, namely those knots $K$ for which the complement of $K$ in $S^3$ admits a unique complete constant curvature hyperbolic metric.\footnote{
All but $32$ of the $1,701,936$ knots up to $16$ crossings are hyperbolic~\cite{hoste1998first}.
The Adams conjecture, which posits that the probability of a prime knot being hyperbolic approaches unity as the crossing number goes to infinity, contradicts another, more widely accepted conjecture that the crossing number of a composite knot is not less than that of each of its factors~\cite{v2020question}.
The relationships one observes sometimes depend crucially on the datasets investigated.
}
Examples of such hyperbolic invariants include the volume of the complement and the lengths of certain longitude and meridian cycles in this metric.
In addition, there are various three-dimensional numerical invariants such as the bridge index and Turaev genus and, since every knot bounds orientable embedded surfaces, invariants that rely on these surfaces such as the determinant, Arf invariant, and the four-dimensional  smooth slice genus.
The physics interpretation of certain knot invariants --- most notably the Khovanov polynomial --- is as well natural in diverse dimensions~\cite{Gaiotto:2011nm, Witten:2011zz}.

Given this vast zoo of knot invariants, we would like to formulate general theories that explicate the relationships between the various quantities.
The generalized volume conjecture, which relates evaluations of the $n$-colored Jones polynomial at roots of unity $J_n(e^{2\pi i/n};K)$ with the hyperbolic volume $\Vol (K)$ and the Chern--Simons invariant $\text{CS}(K)$~\cite{Kashaev1997, Murakami2001, murakami2002kashaev, Gukov:2003na} provides an example of a non-trivial relation and states that
\begin{equation}
    \lim_{n\to\infty} \frac{2\pi \log J_n(e^{2\pi i/n};K)}{n} = \text{Vol}(K) + 2\pi^2 i\, \text{CS}(K) ~,
\label{eq:volume-conjecture}
\end{equation}
where $n$ labels the irreducible representation of $SU(2)$ with dimension $n$.

As an initial step in the endeavor, machine learning supplies a practical tool for identifying how knot invariants are correlated.
For example, neural networks can predicts quasipositivity, the $s$-invariant, and the $\tau$-invariant from braid words together with other input invariants~\cite{hughes2016neural}.
Inspired by the volume conjecture and numerical investigations by Dunfield~\cite{Dunfield2000} and Khovanov~\cite{khovanov2003}, work in~\cite{Jejjala:2019kio} showed that the volume of the knot complement of hyperbolic knots can be machine learned from the Jones polynomial with better than $97\%$ accuracy.\footnote{
The are approximately $840,000$ unique Jones polynomials for knots up to $16$ crossings.
When knots with different volumes have the same Jones polynomial, the volumes differ on average by about $3\%$, so the neural network's performance is nearly optimal on this dataset~\cite{Craven:2020bdz}.
}
This success is explained physically through the analytic continuation of Chern--Simons theory~\cite{Witten:2010cx}.
Indeed, based on an evaluation of the Jones polynomial at $t = e^{3\pi i/4}$, there is a simple formula that approximates the volume to a similar accuracy as the trained neural network~\cite{Craven:2020bdz}.
This constitutes a ``reverse engineering'' of the neural network behavior in~\cite{Jejjala:2019kio}, a task which is known to be difficult in general.
The phase at which the Jones polynomial is evaluated is determined using layer-wise relevance propagation~\cite{montavon2019layer} to ascertain the input feature most important to the neural network's performance.
Likewise, the polynomial invariants also predict the $s$-invariant and the slice genus to $98\%$ accuracy~\cite{Craven:2021ckk}.

The unknot decision problem~\cite{Gukov:2020qaj} and the slice ribbon conjecture~\cite{ghmr} have also recently been attacked using machine learning.
The use of machine learning as a tool to understand the structure of knots more generally has gained popularity~\cite{levitt2019big,kauffman2020rectangular,Hajij:2020alg,pawel2021knot,Hughes:2021ams1,Hughes:2021ams2,vernitski2022reinforcement,grunbaum2022narrowing}.
Perhaps the most impressive result to date in this arena is a theorem that relates the signature $\sigma$, slope, volume, and injectivity radius of hyperbolic knots~\cite{davies2021advancing, davies2021signature}:
\begin{theorem}\label{thm:nature}
There exists a real constant $c$ such that
\begin{equation} 
|2\sigma(K)-\slope(K)| \; \le \; c\, \frac{\Vol(K)}{\inj(K)^3} ~. \label{eq:nature}
\end{equation}
\end{theorem}
\noindent The four topological invariants appearing in Theorem~\ref{thm:nature} were determined from a saliency map that assigned an attribution score to each of a dozen inputs using gradient methods.
Just as the four color theorem~\cite{appel1976every} demonstrated the utility of computers for proving theorems,~\eref{eq:nature} establishes the utility of machine learning for deducing exact results in mathematics.

\comment{
As an aside, we note that results in \cite{davies2021signature} are interesting in that they can also be used to establish an upper bound on $\inj(K)$ in terms of $\Vol(K)$:
\begin{theorem}\label{thm:new}
There exists a real constant $c$ such that
\begin{equation} 
\inj(K)^3 \; \le \; c\, \Vol(K) ~. \label{eq:new}
\end{equation}
\end{theorem}
The constant in Theorem~\ref{thm:new} is that same as the constant in Theorem~\ref{thm:nature}.  Results in \cite{davies2021signature} suggest that $c$ may be as small as $c=0.3$, while the \MH{I'm not actually sure about the constant anymore.  Looking at the proofs closer it looks as though they swap the roles of $c_1$ and $c_2$ in the proof of the main theorem.  I need to confirm the value of the constant.} 
}

In the spirit of these earlier efforts, in this work we systematize the study of how well various numerical knot invariant can be predicted either from a polynomial invariant or from up to three other numerical invariants.
We discover surprising correlations among knot invariants; targeting experiments to search for bipartite and multipartite correlations reveals that the most interesting of our results are actually relations between individual invariants, similar in spirit to the volume conjecture \eqref{eq:volume-conjecture}.
This is rather surprising, as it means that multipartite correlations among knot invariants are either uncommon or subtle compared to relations between single invariants.
The results involve many of the well-known invariants we described earlier, including the Jones polynomial, Floer homological invariants $\tau$ and $\varepsilon$, and a hyperbolic longitude invariant.
These correlations point toward the existence of both known and unknown structures relating algebraic and geometric knot invariants.

Two sections follow.
In Section~\ref{sec:experiments}, we describe our automated exploration of bipartite, tripartite, and more general multipartite correlations within a knot database, including a number of interesting case studies.
In Section~\ref{sec:disc}, we discuss the potential mathematical and physical implications of these case studies.
Two appendices contain information on running experiments and viewing results (Appendix~\ref{app:viewing-results}) and short definitions of selected database invariants (Appendix~\ref{app:invariants}).

\section{Experiments}\label{sec:experiments}

\subsection{Data}
All the experiments were run using the KnotInfo database~\cite{knotinfo}, which we expanded to a larger dataset when possible.
A full list of the 53 invariants which were used in the experiments is included in Appendix~\ref{app:invariants}.
Neural networks were trained to predict given invariants from others, thereby identifying interesting relationships between invariants in the database.
We considered experiments with between one and three inputs to the neural network.
This resulted in over 700,000 possible experiments.
However, not all invariant values are calculated for all knots in the KnotInfo database.
A neural network was only trained if, for a specific combination of invariants, there were more than 1,000 knots to work with. 
The final set of results includes data from around 199,330 experiments.
In some case studies where additional data was available, we supplemented the KnotInfo dataset with all 313,199 hyperbolic knots up to $15$ crossings.

In order to be used in the experiments, some data had to be processed.
All polynomials were simply represented as flattened versions of the vectors appearing in the KnotInfo database.
Any invariants represented by Yes/No or Y/N were converted to ones and zeros, respectively.
For regression tasks, any entries where the output invariant is zero were removed.\footnote{These knots were removed to avoid dividing by zero in the relative mean squared error percentage calculation. Only a handful of knots are removed in these cases, so this removal makes little difference to the experiments.}
For classification tasks, the output invariants were rescaled so that the lowest number in the output data set is zero. 

\subsection{Machine learning architecture}
For each experiment, the neural networks were trained on $80\%$ of the data. The neural networks had three hidden layers, each consisting of one hundred nodes.\footnote{For more detailed discussions of architecture, see~\cite{Jejjala:2019kio,Craven:2021ckk}.}
The ReLu (Rectified Linear Unit) activation function $r(x) \equiv \max(0,x)$ was used on the hidden layers, and the networks were trained using the Adam optimizer. The output layer was either a softmax $\vec{m}(\vec{v})$ with components $m(\vec{v})_k \equiv e^{v_k} / (\sum_j e^{v_j} )$ or a ReLu activation function, depending on whether the task was classification or regression.\footnote{All regression tasks we performed were predicting non-negative quantities, so the ReLu activation could be applied.} 
Similarly, the loss function was either sparse categorical cross-entropy\footnote{The categorical cross-entropy is generally applied after a softmax layer, in which case it is equal to $-\log v_k - \sum_{i\neq k}\log(1-v_i)$ where $v_i$ are the softmax outputs and $k$ represents the true label of the input.  The sparse variant simply allows for the true label dataset to be loaded as integers rather than vectors with a one in the $j^{\text{th}}$ component representing the $j^{\text{th}}$ class.} or mean squared error, depending on the task. Each network was trained for $100$ ``epochs,'' each of which involved a single pass through the entire training set.

\subsection{Case studies}
In this section, we discuss a handful of interesting results. In the cases where the input was a polynomial invariant, in addition to training networks on the full polynomial we also performed a search to find  good evaluations of the polynomial. 
This was done by evaluating the polynomial at points on the square in $\mathbb{C}$ with corners $\pm 1 \pm i$ and training neural networks on each of the evaluations. 
We found that a single evaluation of a polynomial invariant is often sufficient to predict the desired output with high accuracy.\footnote{
For regression tasks, the quantity we refer to as accuracy is obtained from the mean relative error:
$$
\text{accuracy} = 1 - \frac{1}{N} \sum \left| \frac{\text{predicted value}-\text{actual value}}{\text{actual value}} \right| ~,
$$
where $N$ is the size of the dataset.
For classification tasks, it is instead defined simply as $\text{(number of correct classifications)}/\text{(total classifications)}$.  
We modeled discrete invariant prediction as a classification task and continuous invariant prediction as a regression task.
}
In other cases, we found single integer invariants that predicted real-valued quantities, hinting at a very strong underlying correlation.

\subsubsection*{Ozsv\'ath--Szab\'o $\bm{\tau \rightarrow \varepsilon}$}
In an experiment using the expanded dataset, a neural network was able to predict $\varepsilon$ from $\tau$ with essentially perfect accuracy (over $99\%$). This is an expected result, as the sign of $\tau$ is equal to $\varepsilon$ for the majority of these knots. We include this case as a simple example of how the experiments rediscover known results.

\subsubsection*{Jones polynomial $
\bm{\rightarrow \varepsilon}$}

With the exception of a relation between the Jones polynomial and Rasmussen $s$-invariant (which was also uncovered by machine learning techniques \cite{Craven:2021ckk}), the polynomial invariants are not known to relate directly to homological integer invariants like $\varepsilon$ and $\tau$.
In this and the next case study, we find such correlations and discuss their implications further in Section~\ref{sec:disc}.

The Jones polynomial predicted $\varepsilon$ to $96.22\%$ accuracy in the smaller dataset, while the accuracy from predicting the modal value of $\varepsilon$ gives $46.69\%$ accuracy.\footnote{In this and other experiments we compare the accuracy of our neural networks' predictions to a baseline accuracy we obtained by simply predicting the mean (for regression tasks) or mode (for classification tasks) of the target output invariant.}
In the larger dataset, the Jones polynomial predicted $\varepsilon$ to $94.51\%$ accuracy. All experiments that follow in this case study were performed using the KnotInfo dataset.

We are also interested in whether evaluations of the Jones polynomial can be used to predict $\varepsilon$.\footnote{A correlation of this sort, between the Jones polynomial and a homological integer invariant, was discovered for the $s$-invariant in~\cite{Craven:2021ckk}.}
To test this, we chose five random complex points on which to train the neural network.
We then used Layerwise Relevance Propagation (LRP)~\cite{montavon2019layer} to identify relevant inputs. LRP is a technique used to assign a relevance score to each input feature in a neural network, in order to explain how the neural network makes it predictions. It achieves this by modeling relevance as a sort of graph-theoretic flow which propagates backwards through a network, starting at the output layer, and redistributing the flow (relevance scores) into the previous layers while enforcing flow conservation so that the total relevance in each layer is always equal to the original value at the output source.\footnote{See~\cite{Craven:2020bdz} for a more quantitative explanation.} 

In~\cite{Craven:2020bdz}, LRP identified $t = e^{3\pi i/4}$ as a phase at which the evaluation of the Jones polynomial predicts the hyperbolic volume to $97\%$ accuracy.
To pinpoint this phase, the Jones polynomial was evaluated at various roots of unity, in part referencing the volume conjecture and also because of the success of the neural network when given information about the coefficients only and not the degrees.
Here, where there is no \textit{a priori} reason to suspect that phases are important, we consider evaluations at random complex numbers close to the origin in the upper half plane.

\begin{figure}[h!]
    \centering
    \includegraphics[scale=0.5]{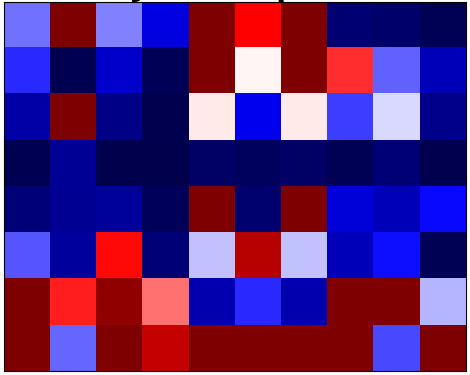}
    \caption{\small{Layerwise Relevance Propagation results from a neural network trained to learn $\varepsilon$ using evaluations of the Jones polynomial. Each column represents a knot, and every pair of rows represent the real and imaginary parts of the Jones polynomial evaluations at a particular point. Red pixels have higher relevance than blue pixels. Notice that certain evaluations are consistently more relevant than others. In this example, the final two rows are the most relevant. Although there are other rows that have high relevance for some knots, at least one of the final two rows (which correspond to a single evaluation) are highly relevant for every knot. These correspond to evaluations of the Jones polynomial at $t = -0.98+0.88i$.}}
    \label{fig:lrp_eps}
\end{figure}

The results of the LRP experiment are shown in Figure~\ref{fig:lrp_eps}. Since one evaluation is consistently more relevant than others, we trained the neural network on this single evaluation of the Jones polynomial. Training on the Jones polynomial at $t = -0.98+0.88i$, the neural network achieved an accuracy of $82.01\%$ at predicting $\varepsilon$ over five runs. 

\begin{figure}[h!]
    \centering
    \includegraphics[scale=0.8]{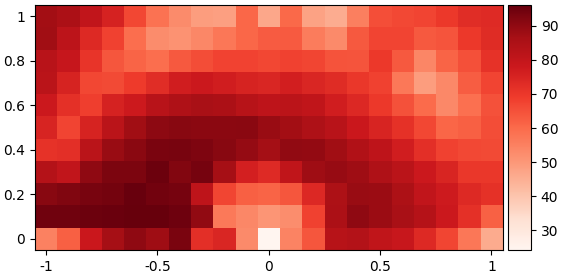}
    \caption{\small{Heatmap showing accuracies of neural network predictions for $\varepsilon$ from single evaluations of the Jones polynomial. The horizontal axis is the real part of the point where the polynomial is evaluated, and the vertical axis is the imaginary part. The colouring shows the accuracy of the neural network predictions. Each point is averaged over five training runs. Only the upper half plane is included, since the Jones polynomial is holomorphic. The best-performing point is at $t = -0.6+0.1i$.}}
    \label{fig:joneseps}
\end{figure}

The LRP experiments indicate that a single evaluation of the Jones polynomial may be sufficient to predict $\varepsilon$. To find the optimal evaluation, we trained the neural network on single Jones polynomial evaluations using points in the square with corners $\pm 1 \pm i$. 
The evaluation of the Jones polynomial at the point $t = -0.6 + 0.1i$ achieved $96.11\%$ accuracy over five runs. Using only the real part of the evaluation, we obtain an accuracy of $87.22\%$. On the larger dataset, the accuracies are $92.92\%$ (full evaluation) and $81.83\%$ (real part of the evaluation). The results of the search for an optimal evaluation are shown in Figure~\ref{fig:joneseps}.

\subsubsection*{Jones/HOMFLY polynomials $\bm{\to}$ Ozsv\'ath--Szab\'o $\bm{\tau}$}
\begin{figure}[h!]
    \centering
    \includegraphics[scale=0.8]{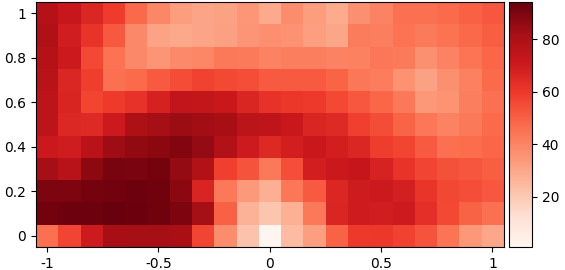}
    \caption{\small{Heatmap showing accuracies of neural network predictions for $\tau$ from single evaluations of the Jones polynomial. The horizontal axis is the real part of the point where the polynomial is evaluated, and the vertical axis is the imaginary part. Each point is averaged over five training runs. Only the upper half plane is included, since the Jones polynomial is holomorphic. The best-performing point is at $t = -0.7+0.1i$.}}
    \label{fig:jonestau}
\end{figure}

The Jones polynomial predicts $\tau$ to $88.86\%$ accuracy, while the accuracy from making predictions from the modal value of $\tau$ on our dataset is $29.21\%$. The HOMFLY polynomial predicts $\tau$ to $88.44\%$ accuracy. Scanning evaluations of the Jones polynomial, we find that the evaluation at $t = -0.7 + 0.1i$ predicts the Ozsv\'ath--Szab\'o $\tau$-invariant to $94.33\%$ accuracy over five runs. This is surprising, since the evaluation at a single number outperforms the full Jones polynomial. Neither the real nor the imaginary parts alone are enough to predict $\tau$ accurately. The magnitude and phase do not give accurate predictions either. The results of the search for an optimal evaluation are shown in Figure~\ref{fig:jonestau}. All experiments for this case study used the KnotInfo dataset.

\subsubsection*{Jones/Kauffman polynomials $\to$ Turaev genus $\bm{g_T}$}

The Jones polynomial is known to be related to the Turaev genus $g_T$ via a bound coming from the span of the polynomial \cite{turaev1990simple}.
More precisely, if $\alpha_{\max}$ and $\alpha_{\min}$ are the maximum and minimum powers of $t$ appearing in the Jones polynomial respectively, the Turaev genus is bounded by
\begin{eqnarray}
    g_T \leq c - |\alpha_{\max} - \alpha_{\min} |\ ,
\label{eq:turaev-bound}
\end{eqnarray}
where $c$ is the crossing number of the knot.
In this case study we find a stronger relationship which uses either the full polynomial or the span in a nontrivial way that avoids using information about the crossing number.

The Jones polynomial predicts $g_T$ to $91.41\%$ accuracy, while the accuracy from predicting the modal value of $g_T$ is $62.71\%$. 
In view of \eqref{eq:turaev-bound}, this is not surprising, since \eqref{eq:turaev-bound} means the Jones polynomial contains some non-trivial information about the Turaev genus.
The amount of information contained by the Jones polynomial is similar to estimating the crossing number itself, though the performance of the network when explicitly given the crossing data is better: we find that training a neural network on $c - |\alpha_{\max} - \alpha_{\min} |$ yields an accuracy of $99.89\%$.\footnote{Simply predicting the upper bound provided by the inequality \eqref{eq:turaev-bound} for the Turaev genus on all knots only yields an accuracy of $65.96\%$.
The improved performance of the neural network given $c - |\alpha_{\max}-\alpha_{\min}|$ may be explained as follows:
on the KnotInfo dataset an accuracy of above $99\%$ can be achieved by guessing $0$ for the Turaev genus if $c-|\alpha_{\max}-\alpha_{\min}| = 0$, and guessing $1$ otherwise.}
Training the neural network on only the span of the Jones polynomial achieves an accuracy of $90.47\%$, which indicates that the span contains nearly as much information as the full Jones polynomial.
One rather trivial explanation for this may be the following: the distribution of spans of the Jones polynomial correlates with crossing number in our dataset.
This would explain the strong performance of the network but would not point to any deep underlying relationship like \eqref{eq:turaev-bound} which allows one to avoid the crossing number.

The Kauffman polynomial predicts $g_T$ to $91.26\%$ accuracy.
As the Jones polynomial is a specialization of the Kauffman polynomial, this is to be expected.
Most evaluations of the Jones polynomial that were tested do not perform better than the baseline $62\%$ accuracy. The best performing evaluation we tested was at the point $t = -1 + 0.2i$, which predicts the Turaev genus to $78.89\%$ accuracy over five runs.

\subsubsection*{Jones/HOMFLY/Conway/Alexander polynomials $\to$ Longitude length $\bm{\ell}$}
With the exception of the generalized volume conjecture, polynomial invariants are not known to be related to hyperbolic invariants.
In this case study we find a very robust correlation between many different polynomial invariants and a particular hyperbolic invariant.
We discuss implications of this result in Section~\ref{sec:disc}.

\begin{figure}
    \centering
    \includegraphics[scale=0.8]{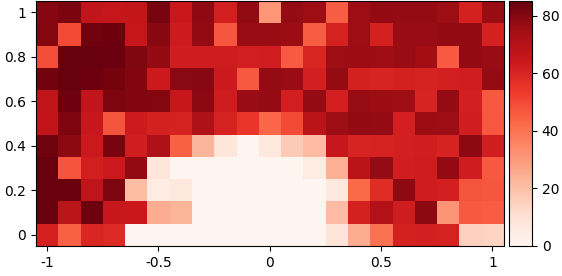}
    \caption{\small{Heatmap showing accuracies of neural network predictions for longitude length from single evaluations of the Jones polynomial. The horizontal axis is the real part of the point where the polynomial is evaluated, and the vertical axis is the imaginary part. Each point is averaged over five training runs. Only the upper half plane is included, since the Jones polynomial is holomorphic. One of the best-performing points is at $t= -1 + 0.2i$.}}
    \label{fig:my_label}
\end{figure}
Using the Jones polynomial, the neural network predicted the longitude length to $89.45\%$ accuracy (averaged over five runs). The accuracy from predicting the mean longitude length for every knot is $46.62\%$, so the Jones polynomial provides a significant improvement.
 
\begin{figure}
    \centering
    \includegraphics[scale=0.45]{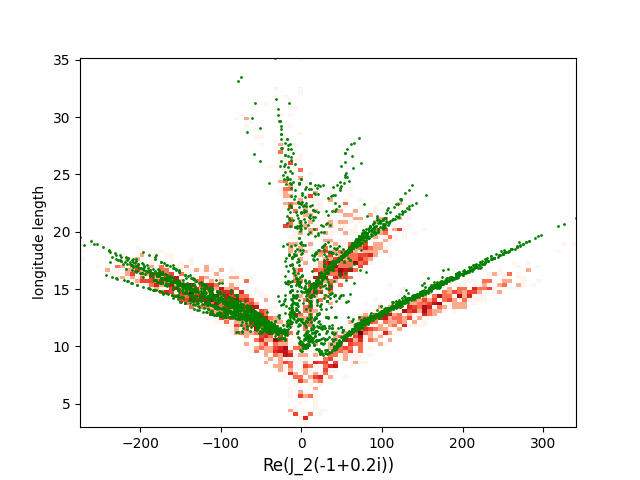}
    \includegraphics[scale=0.45]{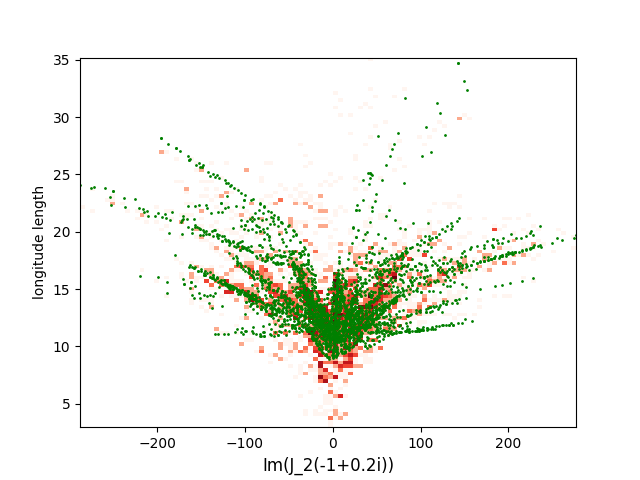}
    \caption{\small{Predictions from the neural network trained on an evaluation of the Jones polynomial at $t = -1 + 0.2i$. \textbf{(Left)} The real part of the Jones polynomial plotted against the longitude length. \textbf{(Right)} The imaginary part of the Jones polynomial plotted against the longitude length. The red points show the actual longitude length (darker red indicating higher density of points) and the green points show the neural network predictions.}}
    \label{fig:nn_preds}
\end{figure}

One can also investigate how well the neural network can perform when training on only an evaluation of the Jones polynomial at a single point. These experiments were performed using the KnotInfo dataset. Many evaluations of those sampled achieve over $80\%$ accuracy. The evaluation of the Jones polynomial at $t = -1 -0.2i$ achieved one of the best accuracies: $87.51\%$ over five runs. On the expanded dataset the evaluation at this point achieved an accuracy of $86.44\%$. We also plot the longitude length against the real and imaginary parts of the Jones polynomial evaluation at $t = -1 - 0.2i$ (Figure~\ref{fig:nn_preds}). Note that the evaluation at $t = -1$ also performs well, signalling a relation with the determinant $\det K$.
The determinant is equal to the absolute value of the Jones polynomial evaluated at $-1$, and can similarily be extracted from the Alexander polynomial.
Since the determinant is a common feature of both polynomials, we might suspect that it is really $\det K$ which serves as a predictor of the longitude length.
As we will see in the next case study, $\det K$ does perform well, but at least a bit worse than the full polynomial invariants.
So the full polynomials contain some extra information about $\ell$ when compared to $\det K$.

\begin{figure}
    \centering
    \includegraphics[scale=0.45]{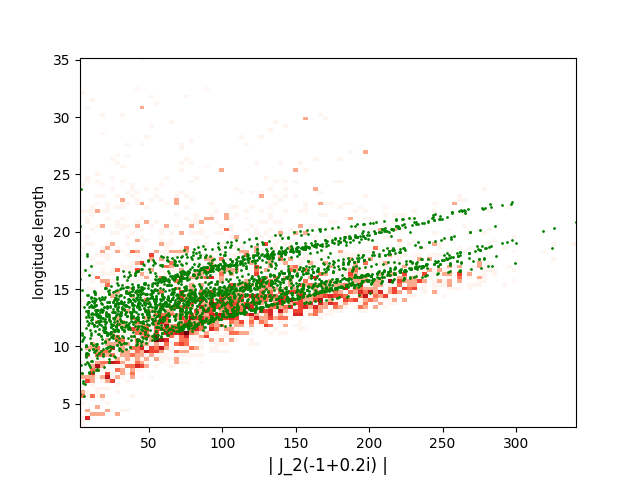}
    \includegraphics[scale=0.45]{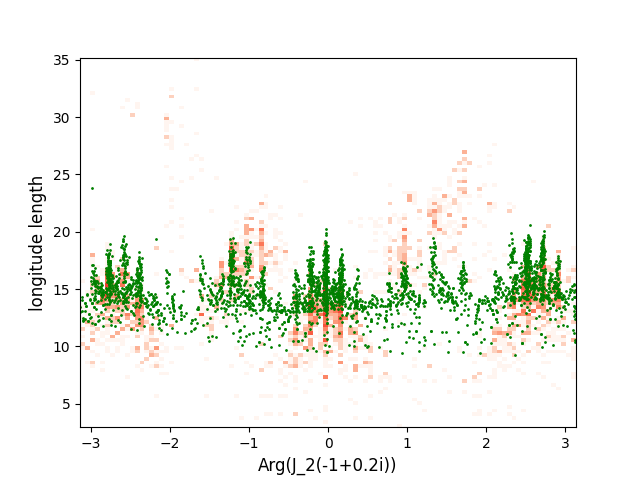}
    \caption{\small{Predictions from the neural network trained on an evaluation of the Jones polynomial at $t = -1 + 0.2i$. \textbf{(Left)} The magnitude of the Jones polynomial plotted against the longitude length. \textbf{(Right)} The phase of the Jones polynomial plotted against the longitude length.  The red points show the actual longitude length (darker red indicating higher density of points) and the green points show the neural network predictions. The plot of the magnitude is very similar in form to the $\det K$ correlation in Figure~\ref{fig:one_ins}.
    }
    }
    \label{fig:nn_preds2}
\end{figure}

We repeated the above experiments using the magnitude and phase of the evaluations as inputs, rather than the real and imaginary parts. The performance dropped to around $82.82\%$, averaged over five runs. Using this form of the evaluation may obscure some information about the longitude length. 
In Figure~\ref{fig:nn_preds2}, we plot the longitude length against the magnitude and phase of the Jones polynomial evaluation at $t = -1-0.2i$. 

The Conway, HOMFLY, and Alexander polynomials all give good predictions of the longitude length as well ($88.73\%$, $88.47\%$, and $87.18\%$ accuracy, respectively).

\subsubsection*{Longitude length predictions}
Initial experiments indicate that certain integer invariants, namely the three genus $g_3$, Ozsv\'ath--Szab\'o $\tau$-invariant, Rasmussen $s$-invariant, and determinant $\det K$, all predict the longitude length $\ell$ quite well. 
These invariants are not known to be related to hyperbolic invariants, as we discuss in Section~\ref{sec:disc}.

\begin{figure}
    \centering
    \includegraphics[scale=0.4]{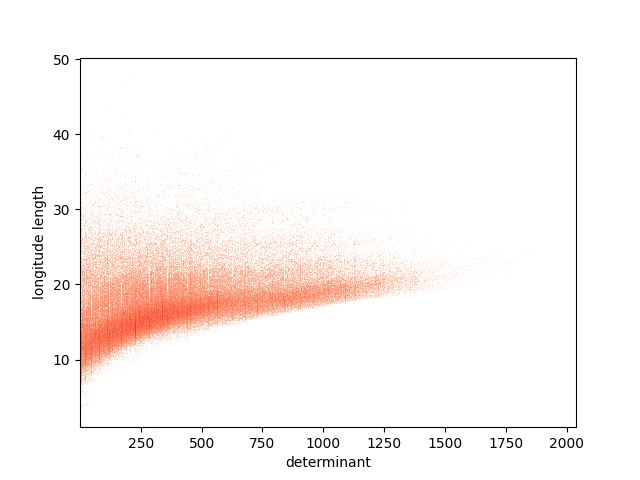}
    \includegraphics[scale=0.4]{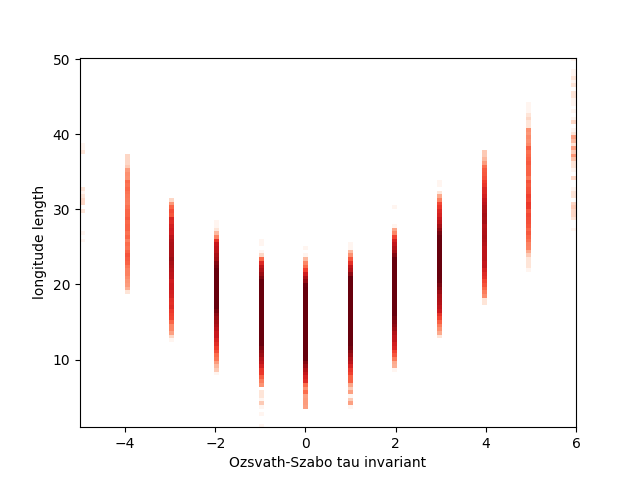}
    \includegraphics[scale=0.4]{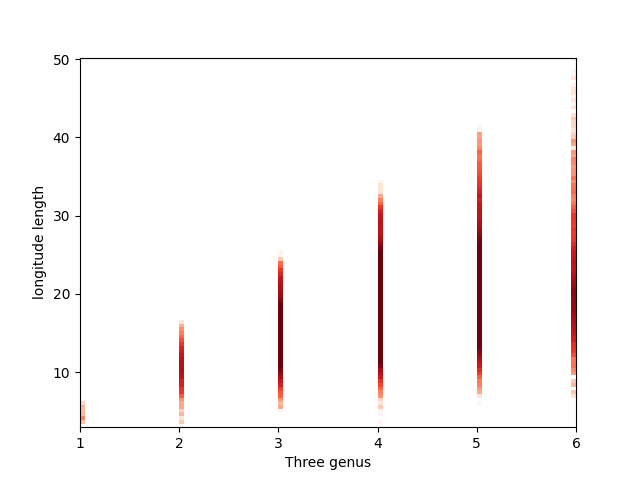}
    \caption{\small{The determinant, Ozsv\'ath--Szab\'o $\tau$-invariant, and the three genus plotted against the longitude length. Darker regions correspond to a higher density of knots.  }}
    \label{fig:one_ins}
\end{figure}

To explore these correspondences, we expand our experiments to include the larger dataset as well as knots with $16$ crossings, for a total of over 1.7 million knots.
Training on only $10\%$ of the data for $100$ epochs, we find that $\det K$, $\tau$, and $g_3$ predict $\ell$ with $85.00\%$, $86.97\%$, and $84.51\%$ accuracy, respectively. Plots of these invariants against the longitude length are shown in Figure~\ref{fig:one_ins}. In Table~\ref{tb:means}, we show the average longitude length for each value of $\tau$ and $g_3$ in the dataset. If we simply predict the mean value of $\ell$ for a particular $\tau$ or $g_3$, we get $86.89\%$ and $84.68\%$ accuracy, respectively.
This recovers the accuracy of the neural network.
Indeed, it explains the majority of the strong performance of the polynomial invariants in predicting $\ell$ as well, since $\det K$ may be extracted from those invariants by a simple evaluation. 

\begin{table}
\centering
\comment{
\begin{tabular}{|c|c|c||c|c|c|}
\hline
    $\tau$ & Mean $\ell$ & $\sigma(\ell)$ & $g_3$ & Mean $\ell$ & $\sigma(\ell)$  \\
    \hline
    -5&32.144&3.45&-&-&-\\
    \hline
    -4&27.54&4.18&-&-&-\\
    \hline
    -3 & 23.08 & 3.37 & - & -& -\\
    \hline
    -2 & 19.31 & 2.77 & - & - &-\\
    \hline
    -1 & 16.33 & 2.61 & - & -& -\\
    \hline
    0 & 15.10& 2.56& - & -- & - \\ 
    \hline
    1&16.55&2.50&1&4.79&0.76\\
    \hline
    2& 19.52 & 2.61 & 2&10.74&1.76\\
    \hline
    3&23.21&3.02&3&14.62&2.62\\
    \hline
    4 & 27.25 & 3.52 & 4 & 17.58&3.43\\
    \hline
    5 & 31.92 & 4.07 & 5 & 19.56& 4.00\\
    \hline
    6 & 37.68 & 5.03 & 6&20.58&4.85\\
    \hline
\end{tabular}}
\begin{tabular}{|c|c||c|c|}
\hline
    $\tau$ & Mean $\ell$ & $g_3$ & Mean $\ell$ \\
    \hline
    $-5$ & $32.14 \pm 3.45$ &$-$&$-$\\
    \hline
    $-4$ & $27.54 \pm 4.18$ &$-$&$-$\\
    \hline
    $-3$ & $23.08 \pm 3.37$ &$-$&$-$\\
    \hline
    $-2$ & $19.31 \pm 2.77$ &$-$&$-$\\
    \hline
    $-1$ & $16.33 \pm 2.61$ &$-$&$-$\\
    \hline
    $0$ & $15.10 \pm 2.56$ &$-$&$-$\\ 
    \hline
    $1$ & $16.55 \pm 2.50$ & $1$ & $4.79 \pm 0.76$\\
    \hline
    $2$ & $19.52 \pm 2.61$ & $2$ & $10.74 \pm 1.76$ \\
    \hline
    $3$ & $23.21 \pm 3.02$ & $3$ & $14.62 \pm 2.62$ \\
    \hline
    $4$ & $27.25 \pm 3.52$ & $4$ & $17.58 \pm 3.43$ \\
    \hline
    $5$ & $31.92 \pm 4.07$ & $5$ & $19.56 \pm 4.00$ \\
    \hline
    $6$ & $37.68 \pm 5.03$ & $6$ & $20.58 \pm 4.85$ \\
    \hline
\end{tabular}
\caption{\small{Average longitude lengths for given values of $\tau$ and $g_3$. We notice that $|\tau|$ may be sufficient to predict $\ell$ to high accuracy.}}
\label{tb:means}
\end{table}

\subsubsection*{Multipartite correlations}

Thus far, all of the case studies we have discussed  actually involved relations between individual invariants.
In our experiments, we also trained neural networks to predict a single invariant from multiple invariants, and we filtered the experiments to isolate truly multipartite correlations.
For instance, we generated a list of experiments with two input invariants which achieved greater than 80\% accuracy where neither of the two inputs alone could predict the output with greater than 80\% accuracy.
Such a list contains tripartite correlation between the two inputs and one output.  We performed similar experiments looking for invariants that could be predicted well from triples of input invariants, likewise filtering out those that could be explained by simpler correlations. 

Perhaps surprisingly, we were not able to find many interesting examples of multipartite correlation.
The few that did appear at the tripartite level with high accuracy were often of the form $(\text{genus},\varepsilon) \to \text{homological invariant}$, where (genus) could be the smooth or topological four genus, three genus, or others and (homological invariant) could be the $s$- or $\tau$-invariants.
These sorts of correlations are expected in the following sense.
The homological invariants are related to the various genus quantities, and indeed this was one of the original motivations for their invention \cite{rasmussen2010khovanov}.
The only subtlety is that a genus is non-negative while the homological invariants can be negative.
The $\varepsilon$-invariant, which is either $\pm 1$ or zero, seems to give just enough additional information to fix the sign, as in the datasets we used it is very often the case that the sign of $\varepsilon$ matches the sign of the homological invariants $\tau$ or $s$.

Another tripartite correlation which is unusual at first sight involves tasks of the form $(\ell , \varepsilon) \to \text{homological invariant}$, where the homological invariant is again $\tau$ or $s$.
This correlation is similar in form to the one we just discussed: the longitude length gives a reasonable estimate of the magnitude of the homological invariant, and $\varepsilon$ provides the sign.
The theoretical underpinning of this correlation is much less clear, and is on the same footing as our experiments showing that $\ell$ may be predicted from $\tau$.
Given that empirical observation, these tripartite correlations are also not surprising, since we know that for knots in our datasets the distribution of $\ell$ is highly correlated with $\tau$.

\section{Discussion}\label{sec:disc}

We have automated a large number of machine learning experiments which search knot datasets for novel relationships between knot invariants.
In addition to recovering earlier relationships we had either only discovered~\cite{Craven:2021ckk,Craven:2022jic} or both discovered~\cite{Jejjala:2019kio} and explained~\cite{Craven:2020bdz}, we have found new relationships between broad classes of knot invariants.
We now discuss the mathematics and physics of these invariant classes in more detail and speculate on possible interpretations of our results.

\subsection{Khovanov and Floer homology}

The first type of novel relation involves two of the most important knot homology theories: Khovanov and Floer homology.
A homology theory is an algebraic structure consisting of a sequence of vector spaces $C_k$ and ``differentials'' $d_k : C_k \to C_{k+1}$ obeying $d_{k+1} \circ d_k = 0$, which allows one to take the quotients $H_k \equiv \text{ker}(d_{k+1}) / \text{im}(d_{k})$.

Khovanov homology~\cite{khovanov2000} is defined by using vector spaces $C_k$ which are built from ``smoothings'' of a knot diagram.
Each crossing in the diagram can be uncrossed in two different ways, and when a choice is made at every crossing we are left with a set of disjoint circles called a smoothing.
After weighting these circles by tensor products of vector spaces and combining them in a certain way, one arrives at Khovanov's $C_k$.
The Khovanov differential is defined by noticing that the set of all smoothings forms a ``cube of resolutions,'' where edges connect two smoothings that differ by a difference in choice at only a single crossing.
Then the $d_k$ are formed from maps which either ``fuse'' the corresponding tensor factors (if the difference in smoothings corresponds to two circles merging into one) or ``split'' a single tensor factor into two (if the difference in smoothings corresponds to a single circle splitting into two).
As we mentioned in Section~\ref{sec:intro}, the graded Euler characteristic of the resulting homology theory $\Kh_k$ is the unnormalized Jones polynomial.
There are also proposals~\cite{Gukov:2004hz,Witten:2011zz,Aganagic:2020olg,Aganagic:2021ubp} for physical theories with Hilbert spaces that are supposed to be isomorphic (in certain sectors) to Khovanov homology.

Knot Floer homology~\cite{Floer:1988,ozsvath2004holomorphic,rasmussen2003floer} can be constructed rather similarly, via a cube of resolutions~\cite{Ozsvath:2009}, but here the resolutions can include singular points where a crossing is made to overlap at a double point.
The graded Euler characteristic of this homology theory $\KFH_k$ is the Alexander polynomial.
The underlying physical theories are related to $\mathcal{N}=2$ supersymmetric gauge theories~\cite{Seiberg:1994aj,Seiberg:1994rs,Witten:1994cg}.
To our knowledge, there are no known physical explanations for why the gauge and string theories relevant for Khovanov homology~\cite{Gukov:2003na,Witten:2011zz} would be related to the variants of $\mathcal{N}=2$ gauge theories in four dimensions~\cite{Witten:1994cg} relevant for Floer homology.

While these two knot homology theories are not obviously related by physics, there are a few mathematical connections regarding the underlying Lie algebras~\cite{Douglas:2014,Tian:2012} and spectral sequences between the homology theories~\cite{Rasmussen:2005,Dunfield:2005si,Ozsvath:2005,dowlin2018spectral,Beliakova:2022proof}.
In particular, the existence of a connection between Khovanov and Floer homology was suspected already in~\cite{Rasmussen:2005}, and is cited as motivation for the development of the $s$-invariant~\cite{rasmussen2010khovanov}.
A detailed understanding of the meaning of spectral sequences between the Hilbert spaces of supersymmetric gauge theories seems like a prerequisite for any physical explanation of the mathematical relationships between the two homology theories, and preliminary steps toward this goal have been explored~\cite{Gukov:2015gmm}.
However, we emphasize that the known mathematical relations also do not provide an explanation for the connection we have found between the Jones polynomial and Floer homology invariants like $\tau$ or $\varepsilon$.

Such specific relationships between Khovanov and Floer homologies may be different facets of the same underlying mysterious connection referred to as ``the FK correspondence'' in~\cite{Rasmussen:2005}.
With an eye toward understanding how our results extend outside of our chosen dataset, we note that~\cite{Rasmussen:2005} comments that this correspondence (explicitly between the rank of the homology theories and the $s$- and $\tau$-invariants) is known to fail for sufficiently complicated knots.
As such, it is possible that the predictions of $\tau$ and $\varepsilon$ from the Jones polynomial that we have observed may break down when the complexity of the knot becomes great enough to violate the FK correspondence of~\cite{Rasmussen:2005} on average.
However, if these correlations survive, they imply new robust entries in the dictionary of the FK correspondence.

\subsection{Knot polynomials and hyperbolic invariants}

One of the most well-known and influential conjectures in knot theory is the volume conjecture~\eqref{eq:volume-conjecture}~\cite{Kashaev1997,Murakami2001,murakami2002kashaev,Gukov:2003na}.
From a coarse-grained perspective, the volume conjecture relates a quantum algebraic invariant (the colored Jones polynomial) to a classical geometric invariant (the hyperbolic volume of the knot complement).
There is a detailed physical understanding of the mechanism which could give rise to something like the volume conjecture: it involves standard notions of analytic continuation and Picard--Lefschetz theory, but applied to a Feynman path integral instead of an ordinary integral~\cite{Witten:2010cx,Witten:2010zr}.
Indeed, this physical underpinning was used in~\cite{Craven:2020bdz} to give a quantitative explanation of the performance of a neural network~\cite{Jejjala:2019kio} which predicted the volume from the Jones polynomial.
Importantly, the classical volume invariant arises in the physics story as a saddle point value of the action evaluated on a certain ``geometric'' gauge field configuration.

Here, we have found novel connections between quantum algebraic invariants and a hyperbolic invariant which is not known to appear as the saddle point value of some action integral: the hyperbolic longitude length $\ell$.
The algebraic invariants are the standard knot polynomials: Jones, HOMFLY-PT, Alexander, Conway, and Kauffman.
It is curious that this relation with the longitude seems to persist between polynomials which serve as Euler characteristics of different homology theories.
A similar statement cannot be made for the hyperbolic volume.
This may be a manifestation of the same mysterious connection between homology theories discussed above and in~\cite{Rasmussen:2005}.

The most na\"{i}ve way to interpret a strong correlation between the Jones polynomial and the longitude length is as a sort of generalization of the volume conjecture to other aspects of hyperbolic geometry.
On one side of the relation, the necessary quantum algebraic ingredients are the colored Jones polynomials $J_n(t;K)$ with a choice of evaluation point $t_n$ (or possibly multiple such points) and the semiclassical limit $n \to \infty$ with $\lim_{n\to\infty} t_n = t_\infty$ finite.

On the other side of the relation, we cannot simply use the classical action of $SL(2,\mathbb{C})$ Chern--Simons theory, because this already leads to the hyperbolic volume of the knot complement $S^3 \setminus K$.
We need a new ingredient to land on the longitude length instead of the volume, and a natural candidate is the Wilson loop operator.
The usual Wilson loop operator in gauge theories is interpreted roughly as the worldline of a charged particle coupled to the gauge field, but in gravitational theories the analogous object can measure properties of the spacetime like geodesic lengths.
There is a well-known relation between $SL(2,\mathbb{C})$ Chern--Simons theory and gravity~\cite{Witten:1988hc}, so it may not be too much of a stretch to use the $SL(2,\mathbb{C})$ gauge theory Wilson loop operator to compute a geodesic length in the classical limit.

With these building blocks in hand, we can write a very rough conjectural relationship
\begin{eqnarray}
    \lim_{n \to \infty} \left( \log |J_n(e^{2\pi i /n}; K'(K))| - \frac{n}{2\pi}  \text{Vol}(K) \right) \approx \ell \ ,
\end{eqnarray}
where $K'(K)$ is a two-component link formed from $K$ (a Wilson loop in the $SU(2)$ irreducible representation of dimension $n$) and an auxiliary Wilson loop which wraps the longitude of the cusp neighborhood in $S^3 \setminus K$ labeled by the fundamental representation of $SU(2)$.
This auxiliary component should arise from a framing of $K$ with zero self-linking; a self-linked auxiliary component would give rise to a winding geodesic in the cusp neighborhood associated with $K$ in the complement $S^3 \setminus K$.

The intuition for this conjecture is that the knot $K$ has a dimension scaling much faster than the auxiliary component, so its Wilson loop forms a background upon which the auxiliary component can be evaluated.
In the large-$n$ limit, the volume conjecture says that this background is the hyperbolic metric on the complement, so the classical limit of the auxiliary Wilson loop computes the exponential of the longitude length.
This effect is a bit like the backreaction of D-branes which occurs in the AdS/CFT correspondence \cite{Maldacena:1997re}.\footnote{Of course, there are several formal differences between the two situations.  In AdS/CFT, one has two descriptions of the D-branes which are equivalent.  In the volume conjecture, the complex critical point yielding the volume is a required contribution to the analytically continued path integral, and no secondary description which avoids it exists.  The intuition we are using is really more related to the idea of emergent geometry in a semiclassical limit, and this theme appears in both situations.}
In the semiclassical limit, D-branes backreact and source a geometric background upon which strings and other objects such as probe branes propagate.
In this situation, our intuition is that the Wilson loop with a dimension that scales with the Chern--Simons level will form a background, namely the hyperbolic geometry, upon which the lighter Wilson loop will measure a length.
The total classical action is then the sum of the background volume term and the longitude length, so subtracting the volume and removing the classical divergence leaves behind the finite longitude length.

Unfortunately, this conjecture has several serious drawbacks.
The fact that the longitude length is defined by enlarging the cusp neighborhood is not accounted for, but perhaps it emerges naturally as some kind of repulsion between the two loops.
Perhaps more seriously, the proposed formula should imply a similar one for the meridian length, and our experiments have not revealed any such correlation.\footnote{This may be an artifact of our small dataset.  The meridian length is distributed quite tightly compared to the longitude length in our dataset, so there is not much room for improvement in prediction by using a polynomial invariant.  If our dataset had a wider range of meridian lengths, perhaps we would notice an improvement.}
Furthermore, the reason we have included the approximate symbol $\approx$ is because it will generally be difficult to disentangle the $O(n^0)$ terms coming from quantum corrections to the classical value of the action (the hyperbolic volume) from the classical contribution of the longitude length, which would also enter at $O(n^0)$ in $\log |J_n|$.
We do not know of a way to isolate the effect of the auxiliary Wilson loop; this is an interesting issue which warrants further exploration.

\subsection{Integer and hyperbolic invariants}

There is at least one quantity which is known to arise from both Khovanov and Floer theory: the determinant $\det K$.
This invariant is related to evaluations of both the Jones and Alexander polynomials:
\begin{equation}
    \det(K) = |J(-1;K)| = |\Delta(-1;K)| ~.
    \label{eq:detja}
\end{equation}
In light of our discussion of several polynomial invariants predicting the longitude length $\ell$, we might suspect that $\det K$ can also predict $\ell$ successfully.
However, this is not quite true, and there is a good reason that it should not be true: $\det K$ is an integer-valued invariant while $\ell$ is a continuum quantity.
In the volume conjecture, there is a scaling limit which allows a continuum quantity to emerge from a sequence of quantum invariants, but there is no such limit for the determinant.
The polynomial invariants allow a sort of fictitious scaling limit by evaluation at a sequence of points in the complex plane.
In fact, the region of the complex plane which we found was relevant for the Jones polynomial in predicting $\ell$ was quite close to $t = -1$, slightly shifted into the positive imaginary half plane.

Nevertheless, it is surprisingly quite possible to predict $\ell$ using $\det K$.
We found strong performance for the $\ell$ prediction task given $\det K$, $\tau$, or the three genus $g_3$.
This bizarre result suggests yet another mysterious connection, but this time between the integer valued homology invariants and continuum hyperbolic invariants.
Some portion of this correlation may be explained by the relationship we found before between the Jones polynomial and the $s$- or $\tau$-invariants.
However, since the knot polynomials contain a great deal of information that is logically distinct from the more subtle homological invariants like $s$ or $\tau$, we cannot explain the integer invariant connection to hyperbolic invariants with only that relation.

\subsection{Unsupervised learning}

A potentially fruitful direction for future work involves unsupervised learning.
In this work, we automated the supervised learning of a large number of knot invariants.
However, this automation required selection of a small number of invariants on which to train, and the selection of an invariant to predict.
We might instead have set up an unsupervised learning problem where we feed the entirety of the dataset into the machine and ask for interesting intrinsic correlations like cluster formation.
This may require considerably more computational power as well as an extended dataset to extract more robust correlations.
We leave this to future work.

\section*{Acknowledgments}
We thank Jim Halverson, Jen Hom, Allison Moore,  and Fabian Ruehle for discussions.
JC acknowledges support from the South African Research Chairs Initiative of the Department of Science and Innovation and the National Research Foundation and the Rosenbaum-Faber Family Graduate Fellowship.
MH is supported by a grant from the National Science Foundation (LEAPS-MPS-2213295).
VJ is supported by the South African Research Chairs Initiative of the Department of Science and Innovation and the National Research Foundation.
AK is supported by the Simons Foundation through the It from Qubit Collaboration.

\appendix

\section{Running your own experiments and viewing results}\label{app:viewing-results}

The file ``\texttt{learn.py}'' in the \href{https://github.com/JessRachel97/Knot_searcher_experiments}{linked repository} \cite{github} can be used to run your own experiments. It uses the KnotInfo dataset and the same neural network architecture that was used in the experiments described in this paper. To run an experiment, run the command
\texttt{python learn.py <num inputs> <inputs> <output>}. For instance, to learn the hyperbolic volume from Jones polynomials, run the command \texttt{learn.py 1 volume jones\_polynomial\_vector}.

The program \texttt{draw\_output.py}, along with the results text files, can be used to produce \LaTeX \, tables of desired results.
Download the program and the results files and run the command \texttt{python draw\_output.py}.
The program will prompt you to answer questions about which results you want to view.
The example below will produce a \texttt{.tex} and a \texttt{.pdf} file containing Table~\ref{extable}.

The results files, in plain text format, is included in the Github repository. The results for all the experiments can be found in \texttt{results\_\{*\}\_fin.txt}, where \{*\} should be replaced by ``one'', ``two'', or ``three''. 
The \texttt{\{*\}\_pruned.txt} files contain results which have been pruned to ensure true correlation between the full input and the output. In other words, a two-input experiment $A + B \rightarrow C$ will only survive pruning if its performance is improved when compared to the one-input experiments $A \rightarrow C$ and $B \rightarrow C$. The file \texttt{results\_x\_high.txt} are those pruned results whose accuracy is over $90\%$. Due to the small size of the dataset, and the fact that each experiment included in the results was only run once, you may find that the results of your own experiments differ from those in the results files. Training occasionally fails resulting in a low accuracy, when in reality the accuracy is high for the majority of runs.

\begin{table}
\begin{center}
\begin{tabular}{|c|c|c|c|c|}
\hline
Input 1 &Output & Accuracy & Mean/Mode & Number \\\hline
\hline
conway polynomial vector&volume&0.8998461019699341&0.7692073351767865&2970\\
\hline
jones polynomial vector&volume&0.957915108777278&0.7692073351767865&2970\\
\hline
kauffman polynomial vector&volume&0.8843623356506978&0.7692073351767865&2970\\
\hline
homfly polynomial vector&volume&0.8621613622952917&0.7692073351767865&2970\\
\hline
alexander polynomial vector&volume&0.918092352173835&0.7692073351767865&2970\\
\hline
\end{tabular}\end{center}\small{\caption{Example of output from the \texttt{draw\_output} program. \label{extable}}}\end{table}

{\small
\begin{lstlisting}[style=base]
user@path# python draw_output.py
Output file name?
my_test 
How many inputs per experiment? !Experiments can have either one, two, or three inputs!
1 
Minimum accuracy?  !From 0 to 1, minimum accuracy achieved by neural network.!
0.5
Maximum accuracy? 
1
Minimum number of knots? !Minimum number of knots used in experiment. There are 2978 knots in total. Experiments with fewer that 100 knots are not run.!
100
Maximum number of knots? 
3000
Comma separated list of inputs (i.e. volume,alternating,quasipositive) or type [ALL] or  [POLY]
POLY
Comma separated list of outputs or type [ALL]
volume

\end{lstlisting}}

\section{Invariants}\label{app:invariants}

As this paper's intended audience includes non-experts in knot theory, we briefly define the invariants used in our case studies.

\paragraph{Alexander polynomial:}
Historically the first knot polynomial, the Alexander polynomial is defined as
\begin{equation}
    \Delta(t;K) = \det ( V - t V^T ) ~,
\end{equation}
where $V$ is a Seifert matrix of the knot.
The Alexander polynomial enjoys the property that $\Delta(t;K)=\Delta(t^{-1};K)$.
Half the degree of the Alexander polynomial supplies a lower bound for the three genus, which is the minimal genus of a Seifert surface for a knot.
(See the description of the three genus below for more details.)
Knot Floer homology~\cite{ozsvath2004holomorphic} categorifies the Alexander polynomial.

\paragraph{Conway polynomial:}
The Conway polynomial $\nabla(t;K)$ is a reparametrization of the Alexander polynomial:
\begin{equation}
    \Delta(t^2;K) = \nabla(t-t^{-1};K) ~,
\end{equation}
where this relation holds up to multiplication by an overall power of $t$.
In terms of a Seifert matrix $V$ of the knot $K$,
\begin{equation}
    \nabla(t;K) = \det ( t^{1/2} V - t^{-1/2} V^T ) ~.
\end{equation}

\paragraph{Determinant:}
The determinant of a knot is a three dimensional numerical invariant obtained from the Seifert matrix:
\begin{equation}
    \det(K) =| \det ( V + V^T )| ~.
\end{equation}
As stated in~\eref{eq:detja}, the determinant can also be computed from the absolute value of evaluations of the Jones or Alexander polynomials at $-1$.

\paragraph{Epsilon:}
As defined by Hom~\cite{hom2014bordered,hom2014knot}, $\varepsilon(K)$ is an invariant of knots derived from 
comparing a pair of maps $F_\tau$ and $G_\tau$ on certain subquotient  complexes of the knot Floer chain complex $CFK^\infty(K)$.  The maps induced on homology by $F_\tau$ and $G_\tau$ cannot both be trivial, and $\varepsilon (K)$ takes values in $\{-1,0,1\}$ depending on which (or possibly both) of these two induced maps vanish.

If $\varepsilon(K)$ vanishes, then so does $\tau (K)$.  Furthermore, for various classes of knots (\textit{e.g.}, Heegaard--Floer thin knots) the value of $\varepsilon (K)$ is equal to the sign of $\tau (K)$.

\begin{figure}[h]
\begin{center}
\includegraphics[width=0.7\textwidth]{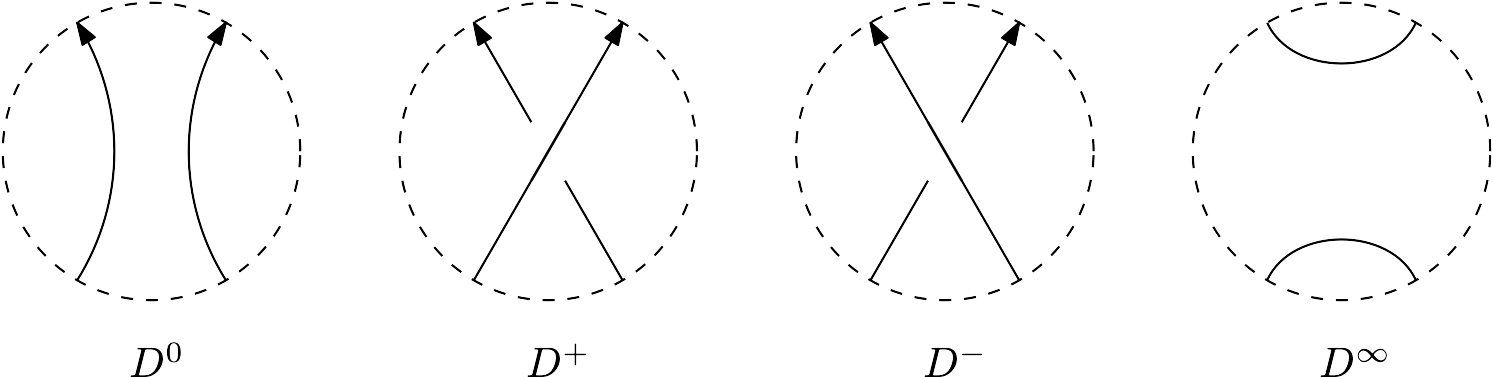}
\end{center}
\caption{Local crossing configurations in $D^0$, $D^+$, $D^-$, and $D^\infty$.}
\label{fig:jonescrossing}
\end{figure}

\paragraph{HOMFLY-PT polynomial:}
Consider the crossing and smoothings on a local region of a knot diagram as in Figure~\ref{fig:jonescrossing}.
Denoting the unknot by $\bigcircle$, the HOMFLY-PT polynomial~\cite{Freyd:1985dx,przytycki2016invariants} is determined by the skein relations:
\begin{itemize}
    \item[I.] $\alpha^{-1} P(\alpha,z;D^+) - \alpha P(\alpha,z;D^-) = z P(\alpha,z;D^0)$~;
    \item[II.] $P(\alpha,z;\bigcircle) = 1$~.
\end{itemize}
The Jones and Alexander polynomials are specializations of the HOMFLY-PT polynomial: 
\begin{eqnarray}
    J(t;K) &=& P(\alpha=t, z=t^{1/2} - t^{-1/2};K) ~, \\
    \Delta(t;K) &=& P(\alpha=1, z=t^{1/2} - t^{-1/2};K) ~,
\end{eqnarray}
where these relations hold up to an overall power of $t$. Analagous to the Jones polynomial discussed below, by taking $\alpha = t^N$, the HOMFLY-PT polynomial has an interpretation in $SU(N)$ Chern--Simons theory~\cite{Witten:1988hf,RamaDevi:1992np}.

\paragraph{Jones polynomial:} The Jones polynomial has been defined in a number of equivalent ways, since its original definition by Jones in terms of von Neumann algebras.  Here, we present a simple recursive definition due to Kauffman~\cite{HKAUFFMAN1987395}.

Let $D^0$, $D^+$, and $D^-$ be three oriented diagrams which are identical except near a single crossing, where they are configured as shown in Figure~\ref{fig:jonescrossing}.  Let $\bigcircle$ denote a knot diagram with no crossings, and let $D \sqcup \bigcircle$ denote a diagram consisting of the union of the diagram $\bigcircle$ and a disjoint knot diagram $D$.  

The Jones polynomial then is a polynomial $J(t;D)$ in the variable $t$ which satisfies the following properties:
\begin{enumerate}
    \item[I.] $t^{-1}J(t;D^+)-tJ(t;D^-) = (t^{1/2}-t^{-1/2})J(t;D^0)$~;
    \item[II.] $J(t;D \sqcup \bigcircle) = -(t^{-1/2}+t^{1/2})J(t;D)$~;
    \item[III.] $J(t;\bigcircle)=1$~.
\end{enumerate}

Although defined in terms of a diagram $D$ of $K$, the resulting polynomial $J(t;D)$ does not depend on the specific choice of $D$, and hence we may refer to the Jones polynomial $J(t;K)$ of the knot $K$.
Khovanov homology~\cite{khovanov2000} categorifies the Jones polynomial.

Witten~\cite{Witten:1988hf} showed that the colored Jones polynomial is computed in Chern--Simons theory on a three manifold $\mathcal{M}$ as follows:
\begin{equation}
J_n(t;K) = \frac{\int [DA]\ e^{iS_\text{CS}[A]}\ W_n(K)}{\int [DA]\ e^{iS_\text{CS}[A]}\ W_n(\bigcircle)} ~,
\end{equation}
where the Wilson loop operator computes the trace of the holonomy along a curve $\gamma$ in the $n$-dimensional representation of $SU(2)$:
\begin{equation}
    W_n(\gamma) = \text{tr}_n\, \mathcal{P} \exp \big( i \oint_\gamma A \big) ~,
\end{equation}
the Chern--Simons action is written in terms of a $\mathfrak{su}(2)$ valued gauge connection:
\begin{equation}
    S_\text{CS}[A] = \frac{k}{4\pi} \int_\mathcal{M} \text{tr}\, \left( A\wedge dA + \frac23 A\wedge A\wedge A \right) ~, \quad k\in \mathbb{Z}^+ ~,
\end{equation}
and $t=e^{2\pi i/(k+2)}$.
The theory is topological, \textit{i.e.}, independent of the metric on $\mathcal{M}$.
Specializing to $n=2$ recovers the Jones polynomial.

\paragraph{Kauffman polynomial:}

\begin{figure}[h]
\begin{center}
\includegraphics[width=0.4\textwidth]{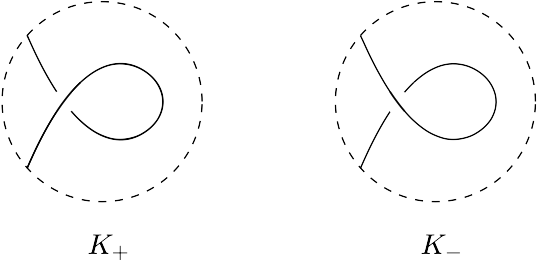}
\end{center}
\caption{Adding a $+1$ kink (left) or a $-1$ kink (right) to obtain the knots $K_+$ and $K_-$, respectively.}
\label{fig:kinks}
\end{figure}

The Kauffman polynomial~\cite{kauffman1990invariant} is a two variable knot polynomial defined as
\begin{equation}
    F(a,z;K)=a^{-w(D)}L(a,z;D)\, .
\end{equation}
Here $w(D)$ is the writhe of a knot diagram $D$, which is the difference between the number of positive crossings and negative crossings, and (referencing Figure~\ref{fig:jonescrossing} with the orientations ignored) $L(a,z;D)$ is determined by skein relations:
\begin{enumerate}
    \item[I.] $L(a,z;D_\pm) = a^{\pm1} L(a,z;D)$~;
    \item[II.] $L(a,z;D^+)+L(a,z;D^-)=z\big(L(a,z;D^0)+L(a,z;D^\infty)\big)$~;
    \item[III.] $L(a,z;\bigcircle)=1$~.
\end{enumerate}
Here, $D_\pm$ adds a $\pm 1$ kink to the diagram $D$ as in Figure~\ref{fig:kinks}.
The Kauffman polynomial is invariant under type 2 and type 3 Reidemeister moves, which, respecitvely, move one loop completely over another and move a string completely over or under a crossing.
The Jones polynomial is a specialization of the Kauffman polynomial:
\begin{equation}
    J(t;K) = F(-t^{3/4},t^{-1/4}+t^{1/4};K) ~.
\end{equation}
The Kauffman polynomial has a definition in $SO(N)$ Chern--Simons theory~\cite{Horne:1989ue}.

\paragraph{Longitude length:}
A cusp of a hyperbolic knot $K$ is a neighborhood of the knot intersected with the knot complement $S^3\setminus K$.  Such a neighborhood can be obtained as the image of a horoball $H$ in the hyperbolic 3-space $\mathbb{H}^3$ under the universal cover $\rho: \mathbb{H}^3 \rightarrow S^3\setminus K$.  For such a horoball $H \subset \mathbb{H}^3$, the preimage $\rho^{-1}(\rho(H))$ will be a family of horoballs, and by increasing the size of $H$, we can arrange for these horoballs to be tangent to each other with disjoint interiors.  In this case the image $\rho (H) \subset S^3\setminus K$ is the maximal cusp of $S^3\setminus K$.

The boundary of a cusp is the torus $T^2$, which has a meridian and a longitude.
The lengths of geodesic representatives of these curves are invariants of $K$.
The maximum longitude length is $5c(K)-6$, where $c(K)$ is the crossing number of the knot~\cite{futer2016survey}.

\paragraph{Ozsv\'ath--Szab\'o $\bm\tau$-invariant:} Any knot $K \subset S^3$ defines a filtration on the Heegaard--Floer chain complex $\widehat{CF}(S^3)$ of $S^3$ (see~\cite{ozsvath2004holomorphic} for the definition of $\widehat{CF}(S^3)$ and~\cite{rasmussen2003floer} or~\cite{ozsvath2003knot} for a definition of the filtration induced by $K$).  If $\mathcal{F}_m(K)$ is the level-$m$ subcomplex induced by the filtration, then the inclusion $\iota_m: \mathcal{F}_m(K) \hookrightarrow \widehat{CF}(S^3)$ induces a map $(\iota_m)_*: H(\mathcal{F}_m(K)) \rightarrow \widehat{HF}(S^3)$ on the homology of these chain complexes.  The Ozsv\'ath--Szab\'o $\tau$-invariant $\tau(K)$ of $K$ is defined to be the minimum integer $m$ for which the map $(\iota_m)_*$ is nontrivial.  The absolute value of $\tau$ gives a lower bound on the slice genus of a knot.

\paragraph{Rasmussen $\bm{s}$-invariant:} In~\cite{lee2008khovanov}, Lee defined a filtered chain complex associated to a knot diagram $D$, by perturbing the differential of the Khovanov chain complex of $D$.  Lee proved that for a knot $K$ the resulting homology always has rank 2, and Rasmussen proved that this homology is always supported in grading $s \pm 1$, where $s = s(K)$ is an invariant of $K$, called the Rasmussen $s$-invariant.  The absolute value of $s$ gives a lower bound on two times the slice genus of a knot.

\paragraph{Three genus:} Every knot $K \subset S^3$ can be expressed as the boundary $K = \partial F$ for some compact, orientable, embedded surface $F \subset S^3$.  Such a surface $F$ is called a Seifert surface for $K$, and the minimal genus $g_3$ among all Seifert surfaces for $K$ is called the three genus (or Seifert genus) of $K$.

Given a Seifert surface $F$ for $K$ of genus $g$, the corresponding Seifert matrix $V$ is a $2g\times 2g$ matrix whose elements are the linking numbers between cycles $a_i$ that are elements  of a basis of $H_1(F)$ and their pushoffs $a_i^+$; \textit{i.e.}, the entries of $V$ are given by $v_{ij} = \text{lk}(a_i,a_j^+)$.
This matrix can be used to define the Alexander and Conway polynomials as well as the determinant of a knot.

\paragraph{Turaev genus:} Given a diagram $D$ of a knot $K$, Turaev~\cite{turaev1990simple} described an algorithm for constructing a closed, orientable, unknotted surface $F(D)\subset S^3$ --- called the Turaev surface of $D$ --- on which  $D$ sits as an alternating diagram.  In particular, if $D$ is an alternating diagram then $F(D)$ is a two-sphere in $S^3$.  The minimal genus of $F(D)$ over all diagrams $D$ of a knot $K$ is called the Turaev genus of $K$, and is denoted by $g_T(K)$.  A knot $K$ is alternating if and only if its Turaev genus is zero.

If $J(t;K)$ is the Jones polynomial of $K$, let $\operatorname{span}(J(t;K))$ denote the difference between the largest and smallest exponents on nonzero monomials in $J(t;K)$, and let $c(K)$ denote the crossing number of $K$ (\textit{i.e.}, the minimal number of crossings in any diagram of $K$).  Then Turaev~\cite{turaev1990simple} proved that the $g_T(K)$ satisfies the following
\begin{equation}
\operatorname{span}(J(t;K)) \leq c(K) - g_T(K) ~.    
\end{equation}

\noindent
\paragraph{Other invariants:}The other 41 invariants we have investigated are tabulated below.
The initial experiments used the data complied in KnotInfo~\cite{knotinfo} for knots up to $12$ crossings.  Invariants for knots up to 16 crossings in the expanded dataset were computed using SnapPy \cite{SnapPy}.
Real invariants are in bold.
Boolean invariants are italicized.
\begin{table}[h!]
{\footnotesize
\begin{tabular}{lll}
arc index &
braid index &
braid length \\
bridge index &
crosscap number &
Morse--Novikov number \\
Nakanishi index &
super bridge index &
Thurston--Bennequin number \\
tunnel number &
unknotting number &
width \\
\textit{Arf invariant} &
$3$d clasp number &
$4$d clasp number \\
smooth $4$d crosscap number &
topological $4$d crosscap number &
smooth concordance crosscap number \\
topological concordance crosscap number &
smooth concordance order &
algebraic concordance order \\
topological concordance order &
smooth concordance genus &
topological concordance genus \\
double slice genus &
smooth four genus $g$ &
topological four genus \\
signature &
\emph{L-space} &
\textbf{Chern--Simons invariant} \\
\textbf{meridian length} &
\textbf{volume} &
\textit{alternating} \\
\textit{fibered} &
\textit{almost alternating} &
\textit{adequate} \\
\textit{quasi-alternating} &
\textit{positive braid} &
\textit{positive} \\
\textit{quasipositive} & 
\textit{strongly quasipositive} &
\end{tabular}
}
\end{table}

\newpage{\small
\bibliographystyle{JHEP}
\bibliography{knots}

\providecommand{\href}[2]{#2}\begingroup\raggedright\begin{thebibliography}{10}

\bibitem{jones85}
V.~F.~R. Jones, \emph{A polynomial invariant for knots via von {N}eumann
  algebras},
  \href{http://dx.doi.org/10.1090/S0273-0979-1985-15304-2}{\emph{Bull. Amer.
  Math. Soc. (N.S.)} {\bf 12} (1985) 103--111}.

\bibitem{HKAUFFMAN1987395}
L.~{H. Kauffman}, \emph{State models and the jones polynomial},
  \href{http://dx.doi.org/https://doi.org/10.1016/0040-9383(87)90009-7}{\emph{Topology}
  {\bf 26} (1987) 395 -- 407}.

\bibitem{Witten:1988hf}
E.~Witten, \emph{{Quantum Field Theory and the Jones Polynomial}},
  \href{http://dx.doi.org/10.1007/BF01217730}{\emph{Commun. Math. Phys.} {\bf
  121} (1989) 351--399}.

\bibitem{Freyd:1985dx}
P.~Freyd, D.~Yetter, J.~Hoste, W.~B.~R. Lickorish, K.~Millett and A.~Ocneanu,
  \emph{{A new polynomial invariant of knots and links}},
  \href{http://dx.doi.org/10.1090/S0273-0979-1985-15361-3}{\emph{Bull. Am.
  Math. Soc.} {\bf 12} (1985) 239--246}.

\bibitem{przytycki2016invariants}
J.~H. Przytycki and P.~Traczyk, \emph{Invariants of links of conway type},
  \href{https://arxiv.org/abs/1610.06679}{{\tt 1610.06679}}.

\bibitem{Alexander:1923}
J.~W. Alexander, \emph{A lemma on systems of knotted curves},
  \href{http://dx.doi.org/10.1073/pnas.9.3.93}{\emph{Proceedings of the
  National Academy of Sciences} {\bf 9} (1923) 93--95},
  [\href{https://arxiv.org/abs/https://www.pnas.org/doi/pdf/10.1073/pnas.9.3.93}{{\tt
  https://www.pnas.org/doi/pdf/10.1073/pnas.9.3.93}}].

\bibitem{alexander1928topological}
J.~W. Alexander, \emph{Topological invariants of knots and links},
  {\emph{Transactions of the American Mathematical Society} {\bf 30} (1928)
  275--306}.

\bibitem{khovanov2000}
M.~Khovanov, \emph{A categorification of the jones polynomial},
  \href{http://dx.doi.org/10.1215/S0012-7094-00-10131-7}{\emph{Duke Math. J.}
  {\bf 101} (02, 2000) 359--426},
  [\href{https://arxiv.org/abs/math/9908171}{{\tt math/9908171}}].

\bibitem{Bar_Natan_2002}
D.~Bar-Natan, \emph{On khovanov’s categorification of the jones polynomial},
  \href{http://dx.doi.org/10.2140/agt.2002.2.337}{\emph{Algebraic \& Geometric
  Topology} {\bf 2} (May, 2002) 337–370},
  [\href{https://arxiv.org/abs/math/0201043}{{\tt math/0201043}}].

\bibitem{ozsvath2004holomorphic}
P.~Ozsv{\'a}th and Z.~Szab{\'o}, \emph{Holomorphic disks and topological
  invariants for closed three-manifolds}, {\emph{Annals of Mathematics} (2004)
  1027--1158}, [\href{https://arxiv.org/abs/math/0101206}{{\tt math/0101206}}].

\bibitem{rasmussen2010khovanov}
J.~Rasmussen, \emph{Khovanov homology and the slice genus}, {\emph{Inventiones
  mathematicae} {\bf 182} (2010) 419--447},
  [\href{https://arxiv.org/abs/math/0402131}{{\tt math/0402131}}].

\bibitem{hom2014bordered}
J.~Hom, \emph{Bordered heegaard floer homology and the tau-invariant of cable
  knots}, {\emph{Journal of Topology} {\bf 7} (2014) 287--326},
  [\href{https://arxiv.org/abs/1202.1463}{{\tt 1202.1463}}].

\bibitem{khovanov2008matrix1}
M.~Khovanov and L.~Rozansky, \emph{Matrix factorizations and link homology},
  {\emph{Fundamenta Mathematicae} {\bf 199} (2008) 1--91},
  [\href{https://arxiv.org/abs/math/0401268}{{\tt math/0401268}}].

\bibitem{khovanov2008matrix2}
M.~Khovanov and L.~Rozansky, \emph{Matrix factorizations and link homology ii},
  {\emph{Geometry \& Topology} {\bf 12} (2008) 1387--1425},
  [\href{https://arxiv.org/abs/math/0505056}{{\tt math/0505056}}].

\bibitem{Gukov:2004hz}
S.~Gukov, A.~S. Schwarz and C.~Vafa, \emph{{Khovanov-Rozansky homology and
  topological strings}},
  \href{http://dx.doi.org/10.1007/s11005-005-0008-8}{\emph{Lett. Math. Phys.}
  {\bf 74} (2005) 53--74}, [\href{https://arxiv.org/abs/hep-th/0412243}{{\tt
  hep-th/0412243}}].

\bibitem{Dunfield:2005si}
N.~M. Dunfield, S.~Gukov and J.~Rasmussen, \emph{{The Superpolynomial for knot
  homologies}},  \href{https://arxiv.org/abs/math/0505662}{{\tt math/0505662}}.

\bibitem{Gukov:2016gkn}
S.~Gukov, P.~Putrov and C.~Vafa, \emph{{Fivebranes and 3-manifold homology}},
  \href{http://dx.doi.org/10.1007/JHEP07(2017)071}{\emph{JHEP} {\bf 07} (2017)
  071}, [\href{https://arxiv.org/abs/1602.05302}{{\tt 1602.05302}}].

\bibitem{Gukov:2017kmk}
S.~Gukov, D.~Pei, P.~Putrov and C.~Vafa, \emph{{BPS spectra and 3-manifold
  invariants}}, \href{http://dx.doi.org/10.1142/S0218216520400039}{\emph{J.
  Knot Theor. Ramifications} {\bf 29} (2020) 2040003},
  [\href{https://arxiv.org/abs/1701.06567}{{\tt 1701.06567}}].

\bibitem{Gukov:2019mnk}
S.~Gukov and C.~Manolescu, \emph{{A two-variable series for knot complements}},
  \href{http://dx.doi.org/10.4171/qt/145}{\emph{Quantum Topol.} {\bf 12} (2021)
  1--109}, [\href{https://arxiv.org/abs/1904.06057}{{\tt 1904.06057}}].

\bibitem{Chandler:2022str}
A.~{Chandler} and E.~{Gorsky}, \emph{{Structures in HOMFLY-PT homology}},
  {\emph{arXiv e-prints} (Sept., 2022) arXiv:2209.13058},
  [\href{https://arxiv.org/abs/2209.13058}{{\tt 2209.13058}}].

\bibitem{hoste1998first}
J.~Hoste, M.~Thistlethwaite and J.~Weeks, \emph{The first 1,701,936 knots},
  {\emph{The Mathematical Intelligencer} {\bf 20} (1998) 33--48}.

\bibitem{v2020question}
A.~V.~Malyutin, \emph{On the question of genericity of hyperbolic knots},
  {\emph{International Mathematics Research Notices} {\bf 2020} (2020)
  7792--7828}, [\href{https://arxiv.org/abs/1612.03368}{{\tt 1612.03368}}].

\bibitem{Gaiotto:2011nm}
D.~Gaiotto and E.~Witten, \emph{{Knot Invariants from Four-Dimensional Gauge
  Theory}}, \href{http://dx.doi.org/10.4310/ATMP.2012.v16.n3.a5}{\emph{Adv.
  Theor. Math. Phys.} {\bf 16} (2012) 935--1086},
  [\href{https://arxiv.org/abs/1106.4789}{{\tt 1106.4789}}].

\bibitem{Witten:2011zz}
E.~Witten, \emph{{Fivebranes and Knots}},
  \href{https://arxiv.org/abs/1101.3216}{{\tt 1101.3216}}.

\bibitem{Kashaev1997}
R.~M. Kashaev, \emph{The hyperbolic volume of knots from quantum dilogarithm},
  \href{http://dx.doi.org/10.1023/a:1007364912784}{\emph{Letters in
  Mathematical Physics} {\bf 39} (1997) 269--275},
  [\href{https://arxiv.org/abs/arXiv:q-alg/9601025}{{\tt
  arXiv:q-alg/9601025}}].

\bibitem{Murakami2001}
H.~Murakami and J.~Murakami, \emph{The colored {J}ones polynomials and the
  simplicial volume of a knot},
  \href{http://dx.doi.org/10.1007/BF02392716}{\emph{Acta Math.} {\bf 186}
  (2001) 85--104}, [\href{https://arxiv.org/abs/arXiv:math/9905075}{{\tt
  arXiv:math/9905075}}].

\bibitem{murakami2002kashaev}
H.~Murakami, J.~Murakami, M.~Okamoto, T.~Takata and Y.~Yokota, \emph{Kashaev's
  conjecture and the chern-simons invariants of knots and links},
  {\emph{Experimental Mathematics} {\bf 11} (2002) 427--435},
  [\href{https://arxiv.org/abs/math/0203119}{{\tt math/0203119}}].

\bibitem{Gukov:2003na}
S.~Gukov, \emph{{Three-dimensional quantum gravity, Chern-Simons theory, and
  the A polynomial}},
  \href{http://dx.doi.org/10.1007/s00220-005-1312-y}{\emph{Commun. Math. Phys.}
  {\bf 255} (2005) 577--627}, [\href{https://arxiv.org/abs/hep-th/0306165}{{\tt
  hep-th/0306165}}].

\bibitem{hughes2016neural}
M.~C. Hughes, \emph{A neural network approach to predicting and computing knot
  invariants}, {\emph{Journal of Knot Theory and Its Ramifications} {\bf 29}
  (2020) 2050005}, [\href{https://arxiv.org/abs/1610.05744}{{\tt 1610.05744}}].

\bibitem{Dunfield2000}
N.~Dunfield, \emph{An interesting relationship between the {J}ones polynomial
  and hyperbolic volume},
  {\emph{\href{https://faculty.math.illinois.edu/~nmd/preprints/misc/dylan/index.html}{Online}}
  (2000) }.

\bibitem{khovanov2003}
M.~Khovanov, \emph{Patterns in knot cohomology, i}, {\emph{Experiment. Math.}
  {\bf 12} (2003) 365--374}, [\href{https://arxiv.org/abs/math/0201306}{{\tt
  math/0201306}}].

\bibitem{Jejjala:2019kio}
V.~Jejjala, A.~Kar and O.~Parrikar, \emph{{Deep Learning the Hyperbolic Volume
  of a Knot}},
  \href{http://dx.doi.org/10.1016/j.physletb.2019.135033}{\emph{Phys. Lett. B}
  {\bf 799} (2019) 135033}, [\href{https://arxiv.org/abs/1902.05547}{{\tt
  1902.05547}}].

\bibitem{Craven:2020bdz}
J.~Craven, V.~Jejjala and A.~Kar, \emph{{Disentangling a Deep Learned Volume
  Formula}}, \href{http://dx.doi.org/10.1007/JHEP06(2021)040}{\emph{JHEP} {\bf
  06} (2021) 040}, [\href{https://arxiv.org/abs/2012.03955}{{\tt 2012.03955}}].

\bibitem{Witten:2010cx}
E.~Witten, \emph{{Analytic Continuation Of Chern-Simons Theory}}, {\emph{AMS/IP
  Stud. Adv. Math.} {\bf 50} (2011) 347--446},
  [\href{https://arxiv.org/abs/1001.2933}{{\tt 1001.2933}}].

\bibitem{montavon2019layer}
G.~Montavon, A.~Binder, S.~Lapuschkin, W.~Samek and K.-R. M{\"u}ller,
  \emph{Layer-wise relevance propagation: an overview},  in \emph{Explainable
  AI: interpreting, explaining and visualizing deep learning}, pp.~193--209.
\newblock Springer, 2019.

\bibitem{Craven:2021ckk}
J.~Craven, M.~Hughes, V.~Jejjala and A.~Kar, \emph{{Learning knot invariants
  across dimensions}},  \href{https://arxiv.org/abs/2112.00016}{{\tt
  2112.00016}}.

\bibitem{Gukov:2020qaj}
S.~Gukov, J.~Halverson, F.~Ruehle and P.~Su\l{}kowski, \emph{{Learning to
  Unknot}},  \href{https://arxiv.org/abs/2010.16263}{{\tt 2010.16263}}.

\bibitem{ghmr}
S.~Gukov, J.~Halverson, C.~Manolescu and F.~Ruehle, \emph{{To appear}},
  \href{https://arxiv.org/abs/22mm.nnnnn}{{\tt 22mm.nnnnn}}.

\bibitem{levitt2019big}
J.~S.~F. Levitt, M.~Hajij and R.~Sazdanovic, \emph{Big data approaches to knot
  theory: Understanding the structure of the jones polynomial},
  \href{https://arxiv.org/abs/1912.10086}{{\tt 1912.10086}}.

\bibitem{kauffman2020rectangular}
L.~Kauffman, N.~Russkikh and I.~Taimanov, \emph{Rectangular knot diagrams
  classification with deep learning},
  \href{https://arxiv.org/abs/2011.03498}{{\tt 2011.03498}}.

\bibitem{Hajij:2020alg}
M.~{Hajij}, G.~{Zamzmi}, M.~{Dawson} and G.~{Muller},
  \emph{{Algebraically-Informed Deep Networks (AIDN): A Deep Learning Approach
  to Represent Algebraic Structures}}, {\emph{arXiv e-prints} (Dec., 2020)
  arXiv:2012.01141}, [\href{https://arxiv.org/abs/2012.01141}{{\tt
  2012.01141}}].

\bibitem{pawel2021knot}
P.~Dłotko, D.~Gurnari and R.~Sazdanovic, \emph{Knot invariants and their
  relations: a topological perspective},
  \href{https://arxiv.org/abs/2109.00831}{{\tt 2109.00831}}.

\bibitem{Hughes:2021ams1}
M.~Hughes, A.~Eubanks and J.~Slone, \emph{Using generative adversarial networks
  to produce knots with specified invariants},  2021.

\bibitem{Hughes:2021ams2}
M.~Hughes, A.~Eubanks and J.~Slone, \emph{Using deep learning to generate knots
  with prescribed invariants},  2021.

\bibitem{vernitski2022reinforcement}
A.~Vernitski, A.~Lisitsa et~al., \emph{Reinforcement learning algorithms for
  the untangling of braids},  in \emph{The International FLAIRS Conference
  Proceedings}, vol.~35, 2022.

\bibitem{grunbaum2022narrowing}
D.~Gr{\"u}nbaum, \emph{Narrowing the gap between combinatorial and hyperbolic
  knot invariants via deep learning}, {\emph{Journal of Knot Theory and Its
  Ramifications} {\bf 31} (2022) 2250003},
  [\href{https://arxiv.org/abs/2204.12885}{{\tt 2204.12885}}].

\bibitem{davies2021advancing}
A.~Davies, P.~Velickovic, L.~Buesing, S.~Blackwell, D.~Zheng, N.~Tomasev
  et~al., \emph{Advancing mathematics by guiding human intuition with ai},
  {\emph{Nature} {\bf 600} (2021) 70--74}.

\bibitem{davies2021signature}
A.~Davies, A.~Juhász, M.~Lackenby and N.~Tomasev, \emph{The signature and cusp
  geometry of hyperbolic knots},  \href{https://arxiv.org/abs/2111.15323}{{\tt
  2111.15323}}.

\bibitem{appel1976every}
K.~Appel and W.~Haken, \emph{Every planar map is four colorable},
  {\emph{Bulletin of the American mathematical Society} {\bf 82} (1976)
  711--712}.

\bibitem{knotinfo}
C.~Livingston and A.~H. Moore, ``Knotinfo: Table of knot invariants.'' URL:
  \url{knotinfo.math.indiana.edu}, November, 2021.

\bibitem{turaev1990simple}
V.~G. Turaev, \emph{A simple proof of the murasugi and kauffman theorems on
  alternating links},  in \emph{New Developments In The Theory Of Knots},
  pp.~602--624.
\newblock World Scientific, 1990.

\bibitem{Craven:2022jic}
J.~Craven, M.~Hughes, V.~Jejjala and A.~Kar, \emph{{(K)not machine learning}},
  in \emph{{Nankai Symposium on Mathematical Dialogues}: {In celebration of
  S.S.Chern's 110th anniversary}}, 1, 2022.
\newblock \href{https://arxiv.org/abs/2201.08846}{{\tt 2201.08846}}.

\bibitem{Aganagic:2020olg}
M.~Aganagic, \emph{{Knot Categorification from Mirror Symmetry, Part I:
  Coherent Sheaves}},  \href{https://arxiv.org/abs/2004.14518}{{\tt
  2004.14518}}.

\bibitem{Aganagic:2021ubp}
M.~Aganagic, \emph{{Knot Categorification from Mirror Symmetry, Part II:
  Lagrangians}},  \href{https://arxiv.org/abs/2105.06039}{{\tt 2105.06039}}.

\bibitem{Floer:1988}
A.~Floer, \emph{{An instanton-invariant for $3$-manifolds}},
  \href{http://dx.doi.org/cmp/1104161987}{\emph{Communications in Mathematical
  Physics} {\bf 118} (1988) 215 -- 240}.

\bibitem{rasmussen2003floer}
J.~A. Rasmussen, \emph{Floer homology and knot complements},
  \href{https://arxiv.org/abs/math/0306378}{{\tt math/0306378}}.

\bibitem{Ozsvath:2009}
P.~Ozsváth and Z.~Szabó, \emph{A cube of resolutions for knot floer
  homology},
  \href{http://dx.doi.org/https://doi.org/10.1112/jtopol/jtp032}{\emph{Journal
  of Topology} {\bf 2} (2009) 865--910},
  [\href{https://arxiv.org/abs/0705.3852}{{\tt 0705.3852}}].

\bibitem{Seiberg:1994aj}
N.~Seiberg and E.~Witten, \emph{{Monopoles, duality and chiral symmetry
  breaking in N=2 supersymmetric QCD}},
  \href{http://dx.doi.org/10.1016/0550-3213(94)90214-3}{\emph{Nucl. Phys. B}
  {\bf 431} (1994) 484--550}, [\href{https://arxiv.org/abs/hep-th/9408099}{{\tt
  hep-th/9408099}}].

\bibitem{Seiberg:1994rs}
N.~Seiberg and E.~Witten, \emph{{Electric - magnetic duality, monopole
  condensation, and confinement in N=2 supersymmetric Yang-Mills theory}},
  \href{http://dx.doi.org/10.1016/0550-3213(94)90124-4}{\emph{Nucl. Phys. B}
  {\bf 426} (1994) 19--52}, [\href{https://arxiv.org/abs/hep-th/9407087}{{\tt
  hep-th/9407087}}].

\bibitem{Witten:1994cg}
E.~Witten, \emph{{Monopoles and four manifolds}},
  \href{http://dx.doi.org/10.4310/MRL.1994.v1.n6.a13}{\emph{Math. Res. Lett.}
  {\bf 1} (1994) 769--796}, [\href{https://arxiv.org/abs/hep-th/9411102}{{\tt
  hep-th/9411102}}].

\bibitem{Douglas:2014}
C.~L. Douglas and C.~Manolescu, \emph{On the algebra of cornered floer
  homology},
  \href{http://dx.doi.org/https://doi.org/10.1112/jtopol/jtt013}{\emph{Journal
  of Topology} {\bf 7} (2014) 1--68},
  [\href{https://arxiv.org/abs/1105.0113}{{\tt 1105.0113}}].

\bibitem{Tian:2012}
Y.~Tian, \emph{A categorification of $u_q$ $sl(1,1)$ as an algebra},
  \href{https://arxiv.org/abs/1210.5680}{{\tt 1210.5680}}.

\bibitem{Rasmussen:2005}
J.~Rasmussen, \emph{Knot polynomials and knot homologies},
  \href{https://arxiv.org/abs/math/0504045}{{\tt math/0504045}}.

\bibitem{Ozsvath:2005}
P.~Ozsv{\'a}th and Z.~Szab{\'o}, \emph{On the heegaard floer homology of
  branched double-covers},
  \href{http://dx.doi.org/10.1016/j.aim.2004.05.008}{\emph{Advances in
  Mathematics} {\bf 194} (June, 2005) 1--33},
  [\href{https://arxiv.org/abs/math/0309170}{{\tt math/0309170}}].

\bibitem{dowlin2018spectral}
N.~Dowlin, \emph{A spectral sequence from khovanov homology to knot floer
  homology},  \href{https://arxiv.org/abs/1811.07848}{{\tt 1811.07848}}.

\bibitem{Beliakova:2022proof}
A.~{Beliakova}, K.~K. {Putyra}, L.-H. {Robert} and E.~{Wagner}, \emph{{A proof
  of Dunfield-Gukov-Rasmussen Conjecture}}, {\emph{arXiv e-prints} (Oct., 2022)
  arXiv:2210.00878}, [\href{https://arxiv.org/abs/2210.00878}{{\tt
  2210.00878}}].

\bibitem{Gukov:2015gmm}
S.~Gukov, S.~Nawata, I.~Saberi, M.~Sto\v{s}i\'c and P.~Su\l{}kowski,
  \emph{{Sequencing BPS Spectra}},
  \href{http://dx.doi.org/10.1007/JHEP03(2016)004}{\emph{JHEP} {\bf 03} (2016)
  004}, [\href{https://arxiv.org/abs/1512.07883}{{\tt 1512.07883}}].

\bibitem{Witten:2010zr}
E.~Witten, \emph{{A New Look At The Path Integral Of Quantum Mechanics}},
  \href{https://arxiv.org/abs/1009.6032}{{\tt 1009.6032}}.

\bibitem{Witten:1988hc}
E.~Witten, \emph{{(2+1)-Dimensional Gravity as an Exactly Soluble System}},
  \href{http://dx.doi.org/10.1016/0550-3213(88)90143-5}{\emph{Nucl. Phys. B}
  {\bf 311} (1988) 46}.

\bibitem{Maldacena:1997re}
J.~M. Maldacena, \emph{{The Large N limit of superconformal field theories and
  supergravity}}, \href{http://dx.doi.org/10.1023/A:1026654312961}{\emph{Adv.
  Theor. Math. Phys.} {\bf 2} (1998) 231--252},
  [\href{https://arxiv.org/abs/hep-th/9711200}{{\tt hep-th/9711200}}].

\bibitem{github}
``Accompanying code.''
  \url{https://github.com/JessRachel97/Knot_searcher_experiments}.

\bibitem{hom2014knot}
J.~Hom, \emph{The knot floer complex and the smooth concordance group},
  {\emph{Commentarii Mathematici Helvetici} {\bf 89} (2014) 537--570},
  [\href{https://arxiv.org/abs/1111.6635}{{\tt 1111.6635}}].

\bibitem{RamaDevi:1992np}
P.~Rama~Devi, T.~R. Govindarajan and R.~K. Kaul, \emph{{Three-dimensional
  Chern-Simons theory as a theory of knots and links. 3. Compact semisimple
  group}}, \href{http://dx.doi.org/10.1016/0550-3213(93)90652-6}{\emph{Nucl.
  Phys. B} {\bf 402} (1993) 548--566},
  [\href{https://arxiv.org/abs/hep-th/9212110}{{\tt hep-th/9212110}}].

\bibitem{kauffman1990invariant}
L.~H. Kauffman, \emph{An invariant of regular isotopy}, {\emph{Transactions of
  the American Mathematical Society} {\bf 318} (1990) 417--471}.

\bibitem{Horne:1989ue}
J.~H. Horne, \emph{{Skein Relations and Wilson Loops in {Chern-Simons} Gauge
  Theory}}, \href{http://dx.doi.org/10.1016/0550-3213(90)90317-7}{\emph{Nucl.
  Phys. B} {\bf 334} (1990) 669--694}.

\bibitem{futer2016survey}
D.~Futer, E.~Kalfagianni and J.~S. Purcell, \emph{A survey of hyperbolic knot
  theory},  in \emph{International Conference on KNOTS}, pp.~1--30, Springer,
  2016.
\newblock \href{https://arxiv.org/abs/1708.07201}{{\tt 1708.07201}}.

\bibitem{ozsvath2003knot}
P.~Ozsv{\'a}th and Z.~Szab{\'o}, \emph{Knot floer homology and the four-ball
  genus}, {\emph{Geometry \& Topology} {\bf 7} (2003) 615--639},
  [\href{https://arxiv.org/abs/math/0301149}{{\tt math/0301149}}].

\bibitem{lee2008khovanov}
E.~S. Lee, \emph{An endomorphism of the {K}hovanov invariant},
  \href{http://dx.doi.org/10.1016/j.aim.2004.10.015}{\emph{Adv. Math.} {\bf
  197} (2005) 554--586}, [\href{https://arxiv.org/abs/math/0210213}{{\tt
  math/0210213}}].

\bibitem{SnapPy}
M.~Culler, N.~M. Dunfield, M.~Goerner and J.~R. Weeks, \emph{Snap{P}y, a
  computer program for studying the geometry and topology of $3$-manifolds}, .

\end{thebibliography}\endgroup
}

\end{document}